\documentclass[11pt,twoside]{article}
\usepackage{atmp}

\usepackage{amsthm}

\theoremstyle{plain}
\newtheorem{theorem}{Theorem}
\newtheorem*{theorem*}{Theorem}
\newtheorem*{mainthm*}{Main Theorem}

\newtheorem{lemma}{Lemma}[section]
\newtheorem*{lemma*}{Lemma}
\newtheorem{proposition}{Proposition}[section]

\newtheorem{corollary}{Corollary}

\theoremstyle{remark}
\newtheorem{remark}{Remark}[section]
\newtheorem*{acknowledgement}{Acknowledgments}

\theoremstyle{definition}
\newtheorem{definition}{Definition}[section]
\newtheorem{example}[remark]{Example}

\numberwithin{equation}{section} 
\usepackage{miscmath} 

\usepackage[ps,dvips,arc,all,knot,frame,curve,poly]{xy}
\usepackage{rg} 
\def\xyc#1\endxyc{{\xy*!C\xybox{#1}\endxy}}
\providecommand{\Xy}{\leavevmode%
 \hbox{\kern-.1em X\kern-.3em\lower.4ex\hbox{Y\kern-.15em}}}

\usepackage{graphicx}
\let\origincludegraphics=\includegraphics
\def\includegraphics{\centering\origincludegraphics}

\usepackage{pretty}
\newrefformat{appendix}{Appendix~\ref{#1}}
\newrefformat{prop}{Proposition~\ref{#1}}
\newrefformat{theorem}{Theorem~\ref{#1}}
\newrefformat{lemma}{Lemma~\ref{#1}}
\newrefformat{xmp}{Example~\ref{#1}}
\newrefformat{definition}{Definition~\ref{#1}}
\newrefformat{cor}{Corollary~\ref{#1}}
\newrefformat{rem}{Remark~\ref{#1}}

\usepackage[newitem,newenum]{paralist}

\newcommand{\expval}[1]{\langle\!\langle#1\rangle\!\rangle}

\newcommand{\RG}[1][]{{\operad{R}_{#1}}}  

\newcommand{\SN}[1][\relax]{\mathcal{S}%
  \ifx#1\relax\relax\else\sb{#1}\fi}
\newcommand{\X}{\mathcal{X}} 

\DeclareMathOperator{\In}{In}
\DeclareMathOperator{\Out}{Out}

\newcommand{\M}{\moduli} 
\newcommand{\Mbar}{\modulibar}
\newcommand{\Mcomb}{\modulicomb}

\newcommand{\Hermitian}[1][\relax]{{\mathcal H}\ifx#1\relax\relax\else(#1)\fi}

\newcommand{\Vertices}[1]{#1^{(0)}}
\newcommand{\Edges}[1]{#1^{(1)}}
\newcommand{\Holes}[1]{#1^{(2)}}

\newcommand{\gint}[2][\Lambda,s_*]{%
\Biggl\langle\!\!\!\Biggl\langle#2\Biggr\rangle\!\!\!\Biggr\rangle_{{#1}}}
\newcommand{\dec}{\sptext{dec}}
\newcommand{\res}{\sptext{res}}

\newcommand{\negquad}{\hskip-.5em\relax}
\newcommand{\negqquad}{\hskip-1.5em\relax}

\begin{document}
\copyrightnotice{2003}{7}{527}{578}

\title{Matrix Integrals and Feynman Diagrams in the Kontsevich Model}

\setcounter{page}{527}
\pagestyle{myheadings}

\url{math.AG/0111082}
 \author{Domenico Fiorenza$^1$ and Riccardo Murri$^2$}
 \address{$^1$Dipartimento di Matematica ``Guido Castelnuovo''\\
Universit{\`a} degli Studi di Roma ``la Sapienza'' \\
p.le Aldo Moro 2, 00185 Roma, Italy}
\address{$^2$Scuola Normale Superiore\\
p.za dei Cavalieri, 7 \\ 56127 Pisa, Italy}
\addressemail{fiorenza@mat.uniroma1.it, riccardo.murri@gmx.it}

\markboth{\it MATRIX INTEGRALS AND FEYNMAN DIAGRAMS\ldots}{\it
  D. FIORENZA AND R. MURRI}

\begin{abstract}
  We review some relations occurring between the combinatorial
  intersection theory on the moduli spaces of stable curves and the
  asymptotic behavior of the 't~Hooft-Kontsevich matrix integrals.  In
  particular, we give an alternative proof of the
  Witten-Di~Francesco-Itzykson-Zuber theorem ---which expresses
  derivatives of the partition function of intersection numbers as
  matrix integrals--- using techniques based on diagrammatic calculus
  and combinatorial relations among intersection numbers. These
  techniques extend to a more general interaction potential.
\end{abstract}

%
%
\cutpage

\section{Introduction}
The aim of this paper is to describe some relations occurring between
combinatorial intersection theory on moduli spaces of stable curves
$\overline{\M}_{g,n}$ and the asymptotic expansion of the matrix
integral
\begin{equation}\label{eq:theintegral}
  \int_{\Hermitian[N]}\exp\left\{
    -\sqrt{-1}\sum_{j=0}^\infty(-1/2)^js_j\frac{\tr
      X^{2j+1}}{2j+1}
  \right\}\ud\mu_\Lambda(X),
\end{equation}
where $\Hermitian[N]$ is the space of $N\times N$ Hermitian
matrices, $\Lambda\in\Hermitian[N]$ is diagonal and positive
definite, and $\mu_\Lambda$ is the Gaussian measure defined
by normalization of $\exp(-\frac{1}{2}\tr\Lambda X^2)\ud X$.

The correspondence between these two apparently unrelated theories is
given by ribbon graphs, which appear on the one side as the cells of an
orbifold cellularization of the moduli spaces of curves, and on the
other side as the Feynman diagrams occurring in the asymptotic
expansion of the Gaussian integral \prettyref{eq:theintegral}.  The
idea to relate the intersection theory on $\overline{\M}_{g,n}$ to the
theory of Feynman diagrams is roughly the following: given a
cohomology class $\omega\in H^*(\overline{\M}_{g,n})$, one seeks a set
of Feynman rules such that, for any ribbon graph $\Gamma$, the
integral of $\omega$ on the cell corresponding to $\Gamma$ equals the
amplitude of $\Gamma$ as a Feynman diagram.  In particular, if one
considers monomials of Miller's classes $\psi_i \in H^2(\M_{g,n})$,
then such a set of rules was found by Kontsevich (see
\cite{kontsevich;intersection-theory;1992}): the Feynman diagrams
expansion of
\begin{equation*}
  \label{eq:cubic}
  \int_{\Hermitian[N]}\exp \Bigl\{
  \frac{\sqrt{-1}}{6} \tr X^3
  \Bigr\} \ud\mu_\Lambda(X)
\end{equation*}
also computes the partition function $Z(t_*)$ of intersection numbers
of the $\psi$ classes, namely,
\begin{equation}
  \label{eq:K1}
  Z(t_*) \biggr\rvert_{t_*(\Lambda)}
  \asymp \int_{\Hermitian[N]}\exp \Bigl\{
  \frac{\sqrt{-1}}{6} \tr X^3
  \Bigr\} \ud\mu_\Lambda(X),
\end{equation}
where 
\begin{equation*}
t_k(\Lambda) := -(2k-1)!!  \tr \Lambda^{-2k-1}
\end{equation*}
are the Miwa coordinates on $\Hermitian[N] / U(N)$. More generally,
one can define combinatorial relatives $\M_{m_*;n}$ of $\M_{g,n}$ and
develop an intersection theory on them; in
\cite{kontsevich;intersection-theory;1992} it is shown that
\begin{equation}\label{eq:highpotential}
Z(s_*;t_*)\biggr\rvert_{t_*(\Lambda)} \asymp
\int_{\Hermitian[N]}\exp\left\{
    -\sqrt{-1}\sum_{j=0}^\infty(-1/2)^js_j\frac{\tr
      X^{2j+1}}{2j+1}
  \right\}\ud\mu_\Lambda(X),
\end{equation}
where $Z(s_*, t_*)$ is the partition function of combinatorial
intersection numbers. Equation \eqref{eq:K1} is recovered as a
special case by setting $s_*=(0,1,0,\dots)$.

In a follow-up \cite{witten;kontsevich-model} to Kontsevich' paper,
Witten proposed a conjecture extending the above relation
\prettyref{eq:K1}: derivatives (up to any order) of the partition
function $Z(s_*;t_*)$ with respect to the $t_*$ variables should admit
a matrix integral interpretation, namely, they should correspond to
asymptotic expansions of Gaussian integrals of the form
\begin{equation}\label{eq:theintegral2}
  \int_{\Hermitian[N]} Q(\tr X,\tr X^3,\tr X^5,\dots) \exp \left\{
    \frac{\sqrt{-1}}{6}\tr X^3
  \right\} \ud\mu_\Lambda(X),
\end{equation}
where $Q$ is a polynomial. Witten also checked the first cases of his
conjecture by Feynman diagrams techniques close to Kontsevich' ones. A
proof of the full statement was later given by Di
Francesco-Itzykson-Zuber
\cite{di-francesco-itzykson-zuber;kontsevich-model}. They solved the
problem explicitly, but their argument rests upon checking some
non-trivial algebraic-combinatorial identities, and does not touch the
geometrical aspect of the problem.

Yet, it is clear by equation \prettyref{eq:highpotential} that the
asymptotic expansions of \prettyref{eq:theintegral2} can be seen as
derivatives of $Z(s_*;t_*)$ with respect to the $s_*$ variables,
evaluated at the point $s_*=(0,1,0,0,\dots)$. Therefore, the Di
Francesco-Itzykson-Zuber theorem is equivalent to the existence of a
linear isomorphism $D:\setC[\del/\del t_*]\to \setC[\del/\del s_*]$, $D: P \mapsto
D_P$ such that, for any differential operator $P(\del/\del t_*)$,
\begin{equation}\label{eq:willfollowby}
  P(\del/\del t_*) Z(s_*;t_*) \biggr\rvert_{s_*=(0,1,0,0,\dots)}
  = D_P(\del/\del s_*) Z(s_*;t_*) \biggr\rvert_{s_*=(0,1,0,0,\dots)}.
\end{equation}
Up to our knowledge, this has first been remarked by
Arbarello-Cornalba in \cite{arbarello-cornalba;dfiz}. From a
geometrical point of view, it means that ---in a certain sense---
combinatorial classes on $\overline{\mathcal{M}}_{g,n}$ are Poincar{\'e}
duals to the Miller classes.  For a formalization of this remark see,
in addition to the already cited papers, also the recent preprints by
Igusa \cite{igusa;miller-morita,igusa;kontsevich-cycles} and Mondello
\cite{mondello;tautological}.

We are going to show how, suitably recasting Witten's computations
from \cite{witten;kontsevich-model} in the language of graphical
calculus, one can prove a statement which generalizes the above
equation \prettyref{eq:willfollowby}. Indeed, denote by $\setC\langle \del/\del
t_*\rangle$ the free non-commutative algebra generated by the $\del/\del
t_*$ (acting on the formal power series in $s_*$ and $t_*$ via its
abelianization $\setC[\del/\del t_*]$); then one can prove the following.
\begin{mainthm*}
  There exist an algebra homomorphism
  \begin{equation*}
    D\colon \setC\langle \del/\del t_*\rangle\to
    \setC\langle\!\langle s_*,\del/\del
    s_*\rangle\!\rangle
  \end{equation*}
  with values in a suitable algebra of formal differential operators in
  the variables $s_*$, such that, for any $P\in\setC\langle \del/\del
  t_*\rangle$,
  \begin{equation}\label{eq:implies}
    P(\del/\del t_*) Z(s_*;t_*)
    = D_P(s_*;\del/\del s_*)Z(s_*;t_*).
  \end{equation}
  Moreover,
  \begin{equation}\label{eq:implies2}
    D_{\del/\del t_k}=c_k{s_1}^{2k+1}\del/\del s_k+ \text{lower
      order terms},\qquad c_k\in \setC.
  \end{equation}
\end{mainthm*}

If $s_*^\circ=(s_0^\circ,s_1^\circ,\dots,s_\nu^\circ,0,0,\dots)$,
where the $s_i^\circ$ are complex constants, then evaluating both
sides of \prettyref{eq:implies} at $s_*=s_*^\circ$ and translating the
result into matrix integral terms, we obtain that there exists a
linear map $Q^{s_*^\circ}:\setC[\del / \del{t_*}]\to \setC[s_*^\circ;
\tr X, \tr X^3, \ldots ]$ such that, for $N\gg 0$,
\begin{align*}
  P(\del / \del{t_*})&\int_{\Hermitian[N]}
  \exp\Biggl\{-\U\sum_{j=0}^\nu(-1/2)^j s^\circ_j
  \frac{\tr(X^{2j+1})}{2j+1}\Biggr\}d\mu_\Lambda(X)=
  \\
  &=\int_{\Hermitian[N]}
  Q^{s_*^\circ}_{P}(X)\exp\Biggl\{-\U\sum_{j=0}^\nu
  (-1/2)^j s^\circ_j
  \frac{\tr (X^{2j+1})}{2j+1}\Biggr\}d\mu_\Lambda(X),
\end{align*}
in the sense of asymptotic expansions. Moreover, equation
\prettyref{eq:implies2} implies that, at the point
$s_*^\circ=(0,1,0,0,\dots)$, the map $Q^{s^\circ_*}$ is a vector space
isomorphism, i.e., there exists a vector space isomorphism \(Q:
\setC[\del/\del{t_*}] \to \setC [\tr X, \tr X^3, \break \tr X^5, \dots]\)
such that, for $N \gg 0$,
\begin{multline*}
  P(\partial_{t_*})\int_{\Hermitian[N]}\exp\Biggl\{\frac{\U}{6}
  \tr X^3\Biggr\}d\mu_\Lambda(X) 
  \\
  = \int_{\Hermitian[N]}\!Q_P(X)\exp\Biggl\{\frac{\U}{6}
  \tr X^3\Biggr\}d\mu_\Lambda(X)
\end{multline*}
in the sense of asymptotic expansions, which is precisely the
statement of the Di~Francesco-Itzykson-Zuber (henceforth referred to
as ``DFIZ'').

\subsection*{Plan of the paper}
The paper is organized as follows.

Section~\ref{sec:kontsevich} contains a brief glossary of
intersection theory on moduli spaces of stable curves and its
combinatorial description; moreover the Kontsevich' Main Identity,
relating the intersection numbers to the combinatorics of ribbon
graphs is recalled.  Finally, the partition function $Z(s_*;t_*)$ of
combinatorial intersection numbers is introduced.

In section~\ref{sec:gc}, ribbon graphs are introduced from a different
point of view, namely, as Feynman diagrams appearing in asymptotic
expansions of certain Gaussian integrals.  The approach to Feynman
diagrams theory is through graphical calculus functors; rules of
graphical calculus are recalled and used extensively all through this
paper.

In section~\ref{sec:tHooft}, a cyclic algebra structure, depending on
a positive definite Hermitian matrix $\Lambda$, is introduced on the
space of $N\times N$ complex matrices; we call it the $N$-dimensional
't~Hooft Kontsevich model.  Feynman rules for this algebra reproduce
the combinatorial terms in Kontsevich' Main Identity.  As a
corollary, it is shown how the partition function $Z(s_*;t_*)$ is an
asymptotic expansion for the partition function of the
't~Hooft-Kontsevich model.

In section~\ref{sec:derivatives}, we use an observation by Witten to
relate first order derivatives of the partition function $Z(s_*;t_*)$
to Laurent coefficients of amplitudes (taken in the
$(N+1)$-dimensional 't~Hooft-Kontsevich model) of graphs with a
distinguished hole.

The long section~\ref{sec:proof} contains the proof of the main
result.  In rough details, it goes as follows.  Laurent coefficients
appearing in Witten's formula of Section~\ref{sec:derivatives} are
polynomials in the eigenvalues of $\Lambda$.  They can be expressed as
amplitudes of ribbon graphs in an extended $N$-dimensional
't~Hooft-Kontsevich model, where vertices are allowed to have
polynomial amplitudes.  One can then apply an recursive procedure to
lower the degree of these polynomials.  In the end, a canonical form
for the expectation values of the graphs related to the first order
derivatives of $Z(s_*;t_*)$ is found.  This canonical form is seen to
equal $D(s_*;\del/\del s_*)Z(s_*;t_*)$ for some $D$ in a certain
(non-commutative) algebra $\setC \langle\!\langle s_*;\del/\del s_*
\rangle\!\rangle$ of power series in $s_*$ and $\del / \del s_*$.  As
a corollary the main result of this paper follows, i.e., the existence
of an algebra homomorphism $D\colon\setC\langle\del/\del t_*\rangle
\to \setC\expval{ s_*;\del/\del s_*}$ such that $D(\del/\del
t_*)Z(s_*;t_*)=D_P(s_*;\del/\del s_*)Z(s_*;t_*)$, for any $P\in
\setC\langle\del/\del t_*\rangle$.

Finally, in Section~\ref{sec:final}, various corollaries and examples
of the main result are given; it is shown how a geometrical
interpretation identifies the combinatorial classes on the
moduli spaces of curves with the Poincar{\'e} duals of the
$\psi$ classes and how a matrix integral translation implies the
Di~Francesco-Itzykson-Zuber theorem
\cite{di-francesco-itzykson-zuber;kontsevich-model}.

\section{Intersection numbers on the Moduli Space of Curves}
\label{sec:kontsevich}

Fix integers \(g\geq 0\), \(n\geq 1\) with \(2-2g-n < 0\). Let
\(\M_{g,n}\) be the moduli space of smooth complete curves of genus
\(g\) with \(n\) marked points \(\{x_1,x_2,\dots,x_n\}\), and let
\(\Mbar_{g,n}\) be its Deligne-Mumford compactification
\cite{deligne-mumford}. The moduli space $ \Mbar_{g,n+1}$ is naturally
isomorphic (as a stack) to the universal curve over $\Mbar_{g,n}$;
that is, if
\begin{equation*}
  \pi: \Mbar_{g,n+1} \to \Mbar_{g,n}
\end{equation*}
is the projection map which ``forgets the marking on the point
$x_{n+1}$'', then the fiber $ \pi^{-1}(p)$ at the generic point $p$ of
$ \Mbar_{g,n}$ has a natural structure of a genus $g$ stable curve
$C_p$ with $n$ marked points $ \{x_1(p),x_2(p), \dots,x_n(p) \}$,
lying in the isomorphism class represented by $p$. So we have $n$
canonical sections
\begin{align*}
  x_i: \Mbar_{g,n}& \to \Mbar_{g,n+1}, \\
  p& \mapsto x_i(p) \in C_p.
\end{align*}   
Define line bundles ${\mathcal L}_i$ on $\Mbar_{g,n}$ by
\begin{equation*}
  \left.{\mathcal L}_i \right|_{p}:=T^*_{x_i(p)}C_p,
\end{equation*}
and denote by $ \psi_i$ the Miller classes (\cite{miller, witten;2dgravity, 
kirwan, morita;structure})
\begin{equation*}
  \psi_i:=c_1({\mathcal L}_i) \in H^2(\Mbar_{g,n}, \setC).
\end{equation*}
Finally, denote by $ \langle \tau_{\nu_1} \cdots \tau_{\nu_n}
\rangle_{g,n}$ the intersection number (\cite{witten;2dgravity})
\begin{equation*}
  \langle \tau_{\nu_1} \cdots \tau_{\nu_n} \rangle_{g,n}:=
\int_{\overline{\mathcal
      M}_{g,n}} \psi_1^{\nu_1} \cdots \psi_n^{\nu_n}.
\end{equation*}
The integral on the right hand side makes sense iff \(\psi_1^{\nu_1}
\cdots \psi_n^{\nu_n} \in \Htop({\Mbar_{g,n}})\), i.e., if and only if
\(\nu_1+ \cdots+ \nu_n=3g-3+n\).  Next, define
\begin{equation*}
  \langle \tau_{\nu_1} \cdots \tau_{\nu_n} \rangle:= \sum_g
  \langle \tau_{\nu_1} \cdots \tau_{\nu_n} \rangle_{g,n}.
\end{equation*}
At most \emph{one} contribution in this sum is non-zero, since \(\langle   
\tau_{\nu_1} \cdots \tau_{\nu_s} \rangle_{g,n}\) can be non-null only   
for \(g=1+ \onethird \left((\sum_{i=1}^n \nu_i)-n \right)\).

\subsection{The generating series and the partition function}
 \label{sec:free-energy-and-partition-function}
 
It is convenient to arrange intersection numbers into some formal series.
\begin{definition}\label{dfn:F-and-Z}
  The generating series of intersection numbers (``free energy
  functional'' in physics literature) is the formal series
  \begin{equation}
    \label{eq:dfn-F}
    F(t_*) := \sum_{g, n} 
    \biggl(1/n! \sum_{\nu_1, \dots,\nu_n}\langle\tau_{\nu_1} \cdots
      \tau_{\nu_n}\rangle_{g, n}^{} t_{\nu_1}\cdots
      t_{\nu_n}\biggr),
  \end{equation}
  
  The partition function is the formal series
  \begin{equation}
    \label{eq:dfn-Z}
    Z(t_*):=\exp F(t_*).
  \end{equation}
\end{definition}  
As remarked by Witten \cite{witten;2dgravity}, algebraic relations
among the intersection numbers $\langle\tau_{\nu_1} \cdots \tau_{\nu_n}\rangle$ translate into
differential equations satisfied by the formal series $F(t_*)$ and
$Z(t_*)$.  In fact, Witten conjectured \cite{witten;2dgravity} and
Kontsevich proved \cite{kontsevich;intersection-theory;1992} that
$\del^2 F (t_*)/ \del {t_0}^2$ satisfies the KdV hierarchy. Moreover,
$F(t_*)$ satisfies the \emph{string equation}
\begin{equation*}
\frac{F(t_*)}{\del
t_0}=\frac{{t_0}^2}{2}+\sum_{i=0}^\infty t_{i+1}\frac{F(t_*)}{\del
t_i}; 
\end{equation*}
thus Kontsevich' result is equivalent to saying that $Z(t_*)$ is a
weight $0$ vector for a Virasoro algebra of differential operators. In
this Virasoro algebra formulation, the Kontsevich-Witten result has
been proven by Witten \cite{witten;kontsevich-model}.

A solution of the KdV hierarchy can be recursively computed, so the
Kontsevich-Witten theorem allows one to recursively compute all
intersection indices. Further details on integrable hierarchies
related to intersection theory on the moduli spaces of stable curves
and, more in general, of stable maps, can be found in
\cite{dubrovin-zhang,EHX,getzler;virasoro,kontsevich-manin,
  okounkov-pandharipande}.

\subsubsection{Kontsevich' Matrix Integral}\label{sec:KMatrI}
Both Kontsevich' and Witten's proofs are based on the integral
representation for the partition function $Z(t_*)$, found by
Kontsevich in \cite{kontsevich;intersection-theory;1992}. Let
$\Hermitian[N]$ be the space of $N\times N$ Hermitian matrices, and
let $\Lambda\in\Hermitian[N]$ be a diagonal matrix with positive real
eigenvalues $\{\Lambda_i\}_{i=1, \dots, N}$. Denote by
$\ud\mu_\Lambda$ the Gaussian measure on $\Hermitian[N]$ obtained by
normalizing $\exp\{-\frac{1}{2}\tr\Lambda X^2\}\ud X$, and let
$t_k(\Lambda)=-(2k-1)!!\tr\Lambda^{-(2k+1)}$. Then:
\begin{equation}
  \label{eq:Z-integral}
  Z(t_*)\bigr\rvert_{t_*(\Lambda)} \asymp \int_{\Hermitian[N]} 
  \exp \Bigl\{ \frac{\U}{6}\tr X^3 \Bigr\}
  \ud\mu_\Lambda(X),
\end{equation}
which holds in the sense of asymptotic expansions as eigenvalues of
$\Lambda$ tend to $\infty$.  The proof of this formula relies on the r{\^o}le played
by ribbon graphs on both sides. Indeed, they arise on the left hand
side as the cells of a combinatorial description of the moduli space
of curves, and on the right hand side as Feynman diagrams in the
asymptotic expansion of the integral. We will now briefly recall the
definition of ribbon graphs and the combinatorial cellularization of
the moduli spaces of curves; ribbon graphs as Feynman diagrams will be
described in \prettyref{sec:gc}.

\subsection{A Triangulation of the Moduli Space of Curves}
A well-known construction (see \cite{harer;cohomology-of-moduli,%
  harer-zagier;euler-characteristic,%
  kontsevich;intersection-theory;1992,%
  looijenga;cellular-decomposition,%
  mulase-penkava}), based on results
of Jenkins-Strebel, leads to a combinatorial description of the moduli
space \(\M_{g,n}\). Let us recall its main points.

\begin{definition}
  \label{dfn:rg} 
  A ribbon graph is a $1$-dimensional CW-complex such that any vertex
  is equipped with a cyclic order on the set of incident half-edges.
  Morphisms of ribbon graphs are morphisms of CW-complexes that
  preserve the cyclic ordering at every vertex.
\end{definition}
Isomorphisms of ribbon graphs are, in particular, homeomorphisms of
the underlying CW-complexes. Call $\Aut\Gamma$ the group of
automorphisms of the ribbon graph $\Gamma$.

\begin{figure}[htbp]
  \includegraphics{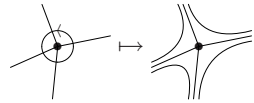}
  \caption{Fattening edges at a vertex with cyclic order.}
  \label{fig:fattening-edges}
\end{figure}

For any ribbon graph \(\Gamma\), denote \(\Vertices{\Gamma}\) the set
of its vertices and \(\Edges{\Gamma}\) the set of its edges.  
Given any ribbon graph, one can use the cyclic order on the vertices
to ``fatten'' edges into thin ribbons\footnote{Hence the name ``ribbon
  graph''.} (see \prettyref{fig:fattening-edges}).  Therefore, a
closed ribbon graph $\Gamma$ is turned into a compact oriented surface
with boundary \(S(\Gamma)\).  The boundary components of $S(\Gamma)$
retract onto particular $1$-homology cycles on $\Gamma$, which we call
``holes''; the set of holes of \(\Gamma\) is denoted $\Holes{\Gamma}$.

The number of boundary components \(n\) and the genus \(g\) of the
closed ribbon graph \(\Gamma\) are defined to be those of the surface
\(S(\Gamma)\).
\begin{definition}\label{dfn:rg-h-l}
  A metric on a closed ribbon graph \(\Gamma\) is a function
  \(\ell:\Edges{\Gamma}\to \setR_{>0}\). A numbering on \(\Gamma\) is
  a map \(h:\Holes{\Gamma}\to \{1,2,\dots,n\}\). 

  Morphisms of numbered (resp.\ metric) ribbon graphs are morphisms of
  ribbon graphs that, in addition, preserve the numbering (resp.\ the
  metric).
\end{definition}
Note that automorphisms of a numbered graph \((\Gamma,h)\)
act trivially on the set \(\Holes{\Gamma}\).

Fix a connected numbered graph \((\Gamma,h)\) with all vertices of valence
$\geq 3$. The set \(\Delta(\Gamma, h)\) of all metrics on $\Gamma$ has a natural
structure of a topological \(m\)-cell (\(m=\card{\Edges{\Gamma}}\))
equipped with a natural action of \(\Aut(\Gamma, h)\); cells \(\Delta(\Gamma, h)\),
with \((\Gamma, h)\) ranging over numbered ribbon graphs of given genus
\(g\) and number of holes \(n\), can be glued to form an
orbi-cell-complex \(\Mcomb_{g,n}\) (see, e.g.,
\cite{looijenga;cellular-decomposition}, \cite{mulase-penkava}), of
which \(\M_{g,n}\) is a deformation retract.
\begin{proposition}[\cite{looijenga;cellular-decomposition}]
  There is an orbifold isomorphism \(\M_{g,n}\times\setR^n_{>0} \simeq
  \Mcomb_{g,n}\).
\end{proposition}
It follows by the above proposition that any integral over
$\M_{g,n}\times\setR^n_{>0}$ can be written as a sum over numbered ribbon
graphs:
\begin{equation}
  \label{eq:cornerstone}
\int_{\Mcomb_{g,n}}
= \sum_{(\Gamma,h)}\frac{1}{\card{\Aut(\Gamma,h)}}\int_{\Delta(\Gamma,h)}
\end{equation}
where $(\Gamma,h)$ ranges over the set of isomorphism classes of
closed connected numbered ribbon graphs of genus $g$ with $n$ holes,
and the factors $1/\card{\Aut(\Gamma,h)}$ appear since we are
integrating over an orbifold \cite{satake;v-manifolds}.  In
\cite{kontsevich;intersection-theory;1992}, representative $2$-forms
\(\omega_i\) for the cohomology classes \(\psi_i\) are found. In terms
of \(\omega_i\)'s, the intersection indices are written
\begin{equation}\label{eq:intersection-indices}
  \langle\tau_{\nu_1} \cdots
  \tau_{\nu_n}\rangle_{g,n}=\int_{\M_{g,n}^{\text{comb}}}\omega_1^{\nu_1}
  \land\cdots
  \land\omega_n^{\nu_n}\land[\setR_{>0}^n],
\end{equation}
where \([\setR_{>0}^n]\) is the fundamental class with compact support
of \(\setR_{>0}^n\).  

As a consequence of \eqref{eq:cornerstone}, Kontsevich finds the
following remarkable identity.
\begin{proposition}[Kontsevich' Main Identity
  {\cite{kontsevich;intersection-theory;1992}}] For all $n$ and $g$,
  the following formula holds:
  \begin{multline*}
    \label{eq:KMI}
    \sum_{\nu_1, \ldots, \nu_n} \langle\tau_{\nu_1} \cdots
    \tau_{\nu_n}\rangle_{g,n} \prod_{i=1}^n \frac{(2\nu_i - 1)!!}
    {\lambda_i^{2\nu_i + 1}}
    \\
    = \sum_{(\Gamma,h)} \frac{1} {\card{\Aut (\Gamma,h)}}\left(
      \frac{1} {2}\right)^{\card{\Vertices{\Gamma}}} \prod_{l\in
      \Edges{\Gamma}} \frac{2} {\lambda_{h(l^+)} + \lambda_{h(l^-)}}
  \end{multline*}
  where: $\lambda_i$ are positive real variables; for any edge $l$ of
  a ribbon graph $\Gamma$, $l^+, l^- \in \Holes{\Gamma}$ denote the
  (not necessarily distinct) holes $l$ belongs to; $(\Gamma,h)$ ranges
  over the set of isomorphism classes of closed connected numbered
  ribbon graphs of genus $g$ with $n$ holes.
\end{proposition}

\subsection{Intersection theory on combinatorial moduli spaces}
 \label{sec:moduli-of-graphs}
 
Kontsevich described a
natural generalization of these constructions.
\begin{definition}\label{dfn:rg-combinatorial-type}
  Let $m_* := (m_0, m_1, \ldots, m_k, \ldots)$ be a sequence of
  non-negative integers such that $m_i \not= 0$ only for a
  \emph{finite} number of indices $i$.  A ribbon graph \(\Gamma\) is
  said to be of \emph{combinatorial type} \(m_*\) if it has exactly
  \(m_i\) vertices of valence \(2i+1\), for \(i\geq 0\), and no
  vertices of even valence.
\end{definition}
One can consider the set of all cells \(\Delta(\Gamma,h)\), where
\(\Gamma\) ranges over closed connected ribbon graphs of a given
combinatorial type \(m_*\), and $h: \Holes{\Gamma} \to \{1, \ldots,
n\}$ is a numbering on the holes of $\Gamma$.  It can be shown that
these cells can be glued together into an orbifold \(\M_{m_*,n}\).  If
\(m_0=0\), then \(\M_{m_*,n}\) is a sub-orbifold of
\(\cup_g\M_{g,n}^{\text{comb}}\) and the support of a homological
cycle.

Equation \prettyref{eq:intersection-indices} can
be generalized to the following definition.
\begin{definition}\label{dfn:comb-intersection-index}
  The combinatorial intersection index \(\langle\tau_{\nu_1} \cdots
  \tau_{\nu_n}\rangle_{m_*, n}\) is defined by
  \begin{equation*}
    \label{eq:combinatorial-intersection-indices}
    \langle\tau_{\nu_1} \cdots \tau_{\nu_n}\rangle_{m_*, n} :=
    \int_{\M_{m_*, n}^{\text{comb}}} 
    \omega_1^{\nu_1} \land\cdots
    \land\omega_n^{\nu_n}\land[\setR_{>0}^n].
  \end{equation*} 
\end{definition}
Since $\bigcup_{m_1} \M_{(0,m_1,0,\ldots);n} = \bigcup_g
\Mcomb_{g,n}$, then one easily computes:
\begin{equation}
\label{eq:top-dim-stratum}
\sum_{m_1}\langle\tau_{\nu_1} \cdots
\tau_{\nu_n}\rangle^{}_{0,m_1,0,0,\dots;n}=\sum_g\langle\tau_{\nu_1}
\cdots
\tau_{\nu_n}\rangle_{g,n} =: \langle\tau_{\nu_1} \cdots \tau_{\nu_n}\rangle.
\end{equation}

Also Kontsevich' Main Identity can be generalized to this
combinatorial context.
\begin{proposition}[Kontsevich' Main Identity
  {\cite{kontsevich;intersection-theory;1992}}] \label{prop:KMI2}
  For any $n$ and any combinatorial type $m_*$, the
following
  formula holds:
  \begin{multline}
    \label{eq:KMI2}
    \sum_{\nu_1, \ldots, \nu_n} s_*^{m_*} \langle\tau_{\nu_1} \cdots \tau_{\nu_n} \rangle_{m_*;n} \prod_{i=1}^n
    \frac{(2\nu_i - 1) !!} {\lambda_i^{2\nu_i + 1}}
    \\
    = \sum_{(\Gamma,h)} \frac{1} {\card{\Aut (\Gamma,h)}} \prod_{j=0}^\infty\left(\frac{s_j}
      {2^j} \right)^{m_j} \prod_{l\in \Edges{\Gamma}} \frac{2} {\lambda_{h(l^+)} +
      \lambda_{h(l^-)}}
  \end{multline}
  where $(\Gamma,h)$ ranges over the set of isomorphism classes of closed
  connected numbered ribbon graphs with $n$ holes and combinatorial type
$m_*$, and the $s_*$ are complex variables.  
\end{proposition}

The free energy and partition function make sense also in this broader
setting:
\begin{equation}\label{eq:dfn-Fst}
  F(s_*;t_*) = \sum_{m_*, n}F_{m_*, n}(s_*;t_*) 
  := \sum_{m_*, n} \biggl( 
  \frac{1}{n!} \!\sum_{\nu_1, \dots,\nu_n}\!\! 
  s_*^{m_*} \langle\tau_{\nu_1} \cdots \tau_{\nu_n}\rangle_{m_*, n}^{} t_{\nu_1}\cdots t_{\nu_n} 
  \biggr),
\end{equation}%
where \(s_*^{m_*}=\prod_{i=0}^\infty s_i^{m_i}\) (this is actually a
\emph{finite} product), and
\begin{equation}\label{eq:dfn-Zst}
  Z(s_*;t_*):=\exp F(s_*;t_*).
\end{equation}%
Correspondingly, one has the integral representation:
\begin{equation}
  \label{eq:Zst-integral}
  Z(s_*; t_*)\biggr\rvert_{t_* = t_*(\Lambda)} \asymp 
  \int_{\Hermitian[N]} \exp \Bigl\{ -\U \sum_j s_j \frac{\tr X^{2j+1}}
  {2j+1} \Bigr\} \ud\mu_\Lambda,
\end{equation}
which will be proved in \prettyref{prop:Z=Z}.  
We recover the usual relations \eqref{eq:dfn-F},
\eqref{eq:dfn-Z} and \eqref{eq:Z-integral} by setting $s_* =
(0,1,0,0,\dots)$ in \eqref{eq:dfn-Fst}, \eqref{eq:dfn-Zst} and
\eqref{eq:Zst-integral},
because $\bigcup_{m_1} \M_{(0,m_1,0,0,\dots);n} =
\bigcup_g \Mcomb_{g,n}$.

\begin{remark}
  The variables \(s_*\) and \(t_*\) are actually of \emph{two}
  different kinds. Indeed, the \(t_*\) variables are free
  indeterminates, whereas the \(s_*\) variables are the structure
  constants of a 1-dimensional cyclic \(A_\infty\) algebra (see
  \cite{kontsevich;feynman}). Since there are no constraints on the
  structure constants of a 1-dimensional cyclic \(A_\infty\)-algebra,
  the \(s_*\)'s are free in the context of this paper. However, when
  dealing with combinatorial classes arising from higher dimensional
  cyclic \(A_\infty\) algebras, the different nature of the two set of
  variables becomes evident.
\end{remark}

\section{Ribbon graphs as Feynman diagrams}
\label{sec:gc}

In \cite{reshetikhin-turaev;ribbon-graphs}, Reshetikhin and Turaev
defined graphical calculus as a functorial correspondence between
certain sets of graphs and morphisms in suitable categories. In one of
its incarnations, this graphical calculus is suitable for working on
ribbon graphs: we follow our treatment \cite{fiorenza-murri;feynman} and
refer the reader also to
\cite{fiorenza;integration-groupoids,oeckl;braided-qft}
for precise statements and proofs.

\begin{definition}
  \label{dfn:rg-with-legs}
  A ribbon graph with $n$ legs is a ribbon graph with a distinguished
  subset of $n$ univalent vertices, called ``endpoints''. An edge
  stemming from one endpoint is called a ``leg''. Edges which are not
  legs, and vertices which are not endpoints are called ``internal''.
   We shall divide internal vertices into two classes (``colors''):
 ``ordinary'' and ``special'' vertices. Morphisms of ribbon
graphs
  with legs map endpoints into endpoints, and preserve
vertex color.
\end{definition} 
The set of isomorphism classes of ribbon graphs with $n$ legs is denoted
by the symbol $\R(n)$; we also set 
\begin{equation*}
  \R:=\bigcup_{n=0}^\infty\R(n). 
\end{equation*} 
An element of $\R(0)$ is called a \emph{closed ribbon graph}.  
Disjoint union gives a map $\R(m) \times \R(m) \to
\R(m+n)$, hence a map $\R \times \R \to \R$.

We will say just ``vertex'' to mean ``internal vertex''. Moreover,  
by abuse of notation, a connected and simply connected
ribbon graph with exactly one internal vertex will be called simply a
``vertex'' (ordinary or special depending on the color of the internal
vertex). The $n$-valent \emph{special} vertex will be denoted by the
symbol ${\sf v}_n^{}\in\R(n)$.

\begin{figure}[htbp]
  \includegraphics{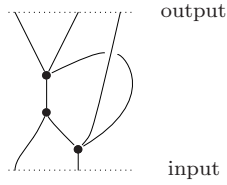}
  \caption{A ribbon graph of type (2,3).}
  \label{fig:ribbon-graph-xmp}
\end{figure}

\begin{definition}
  \label{dfn:rg-with-num-legs}
  A ribbon graph of type $(p,q)$ is a ribbon graph with $p+q$ legs,
  which are partitioned into two disjoint totally ordered subsets: $p$
  ``inputs'' and $q$ ``outputs'' (see
  \prettyref{fig:ribbon-graph-xmp}).
  
  Morphisms of ribbon graphs of type $(p,q)$ are morphisms of ribbon
  graphs with $p+q$ legs which send input legs into input legs, output
  legs into output legs, and preserve the total order on both.
\end{definition}
The set isomorphism classes of ribbon graphs of type $(p,q)$ is
denoted by the symbol $\R(p,q)$ and its $\setC$-linear span by the symbol
$\RG(p,q)$.

The sets of inputs and outputs of a ribbon graph $\hat\Gamma$ of type
$(p,q)$ are denoted, respectively, as $\In(\hat\Gamma)$ and $\Out(\hat\Gamma)$.
The datum of the total order on the legs (of either kind) is
equivalent to a numbering, i.e., to bijections
\begin{equation*}
  \In(\hat\Gamma) \leftrightarrow \{1,\dots,p\},
  \qquad
  \Out(\hat\Gamma) \leftrightarrow \{1,\dots,q\}.
\end{equation*}

One can define a composition product $\R(p,q) \times \R(r,p) \ni (\hat\Phi,
\hat\Psi) \to \hat\Phi \circ \hat\Psi \in \R(r,q)$ by gluing input edges of $\hat\Phi$
with corresponding output edges of $\hat\Psi$. It extends to a bilinear
composition product \(\RG(p,q)\otimes\RG(r,p)\to\RG(r,q)\) which turns \(\RG\)
into the $\Hom$-functor of a category whose objects are natural
numbers; denote this category by the same symbol \(\RG\).
Juxtaposition defines a tensor product $\otimes$ on (isomorphism classes of)
ribbon graphs; one can check that $\otimes$ makes the category $\RG$ monoidal.

\begin{figure}[htbp]
  \includegraphics{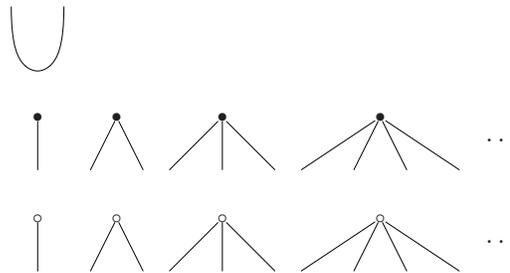}
  \caption{The generators of $\RG$: the bent edge, ordinary
    vertices, special vertices}
  \label{fig:rg-generators}
\end{figure}

Every ribbon graph can be assembled from elementary graphs: vertices
(both ordinary and special) with $k$ inputs and no output for each $k
\geq 1$, and connecting ``bent'' edges with both legs outgoing (see
\prettyref{fig:rg-generators}).  As a monoidal category, $\RG$ is
generated by elements corresponding to these pieces; therefore, one
can define a monoidal functor just by assigning its values on
generating elements.

\begin{definition}
  A cyclic algebra $A = (V, b, T_1, T_2, \ldots)$ over $\setC$ is the
  data of a $\setC$-linear space $V$, of a symmetric non-degenerate
  bilinear form $b: V \otimes V \to \setC$, and of cyclically invariant
  tensors $T_r: V\tp{r} \to \setC$:
  \begin{equation*}
    T_r (v_1 \otimes \ldots \otimes v_{r-1} \otimes v_r) 
    = T_r (v_r \otimes v_1 \otimes \ldots v_{r-1}).
  \end{equation*}
\end{definition}

Let $\langle V \rangle_{x_*}$ be the category having the tensor powers $V\tp{r}$,
for $r \geq 0$, as objects, and
\begin{equation*}
  \Hom_{x_*} (V\tp{p}, V\tp{q}) := \Hom (V\tp{p}, V\tp{q}) \otimes \setC[x_*]
\end{equation*}
as $\Hom$-spaces.  Since $b$ is non-degenerate, it induces a canonical
isomorphism between \(V\) and its dual, and we can give $\langle V
\rangle_{x_*}$ a structure of a rigid monoidal category in which every
object is self-dual.
\begin{proposition}
  \label{prop:gc1}
  Given a cyclic algebra $A$, there is a unique monoidal functor $Z_{A}:
  \RG \to \langle V \rangle_{x_*}$ that maps: \begin{itemize} 
  \item $r$-valent ordinary vertices to
    morphisms $x_rT_r$; 
  \item $r$-valent special vertices to morphisms $T_r$;
  \item bent edges to the copairing \(b\spcheck:
    \setC\to V\otimes V\) dual to the pairing $b: V
    \otimes V \to \setC$.
  \end{itemize}

  The graphical calculus functor $Z_A$ defines (a family of) linear
  maps
  \begin{equation*}
    Z_{A}: \RG(p,q) \to \Hom(V\tp{p}, V\tp{q})\otimes\setC[x_*]
  \end{equation*}
  such that 
  \begin{gather*}
    Z_{A} (\hat\Phi\circ\hat\Gamma) = Z_{A}(\hat\Phi)\circ Z_{A}(\hat\Gamma),
    \\
    Z_{A}(\hat\Phi\otimes\hat\Gamma) = Z_{A}(\hat\Phi) \otimes Z_{A}(\hat\Gamma).
  \end{gather*}
\end{proposition}
Note that, if $\hat\Gamma$ is a ribbon graph of type $(0,0)$, then
$Z_{A}(\hat\Gamma)$ is an element of $\Hom(V\tp{0}, V\tp{0})\otimes\setC[x_*]\simeq
\setC[x_*]$, i.e., it is actually a polynomial with complex coefficients.

\begin{remark}
The vector space $V$ is called the space of ``fields''.  The tensor
$Z_{A}(\hat\Gamma)$ is the ``amplitude'' of the graph $\hat\Gamma$; in the
graphical notation, structure constants of this tensor are denoted by
the graph with indices attached to the legs, whereas the same graph
with no indices will stand for the amplitude tensor itself.
Amplitudes of vertices and bent edges are called, respectively,
``interactions'' and ``propagators''. The data of propagators and
interactions are called the \emph{Feynman rules} of $Z_A$.
\end{remark}

\subsection{Expectation values of graphs}
\label{sec:expectation}

The usual correspondence between Feynman diagrams and Gaussian
integrals will play a key r{\^o}le in this paper; for our purposes,
we can summarize it in the following.

Let $(V_\setR,b)$ be a real Hilbert space and $(V_\setC,b)$ its
complexification; let $A := (V_\setC, b, T_1, T_2, \ldots)$ be a fixed cyclic
algebra structure on $V_\setC$.

For any ribbon graph \(\Psi\in\R\) (possibly with special vertices), we
denote by \(\R_{\Psi}^{}\) the set of (isomorphism classes of) ribbon
graphs containing \(\Psi\) as a distinguished sub-graph and having no
special vertex outside \(\Psi\). By saying that the sub-graph \(\Psi\) is
distinguished, we require that any automorphism of an object \(\Gamma \in
\R_{\Psi}^{}\) maps \(\Psi\) onto itself.  It follows from the definition
that \(\R_\emptyset^{}\) is the set of ribbon graphs having only ordinary
vertices.

The amplitude of an element of $\R(n)$ is not a well-defined tensor,
since there is no distinction between inputs and outputs and no
ordering on the legs. On the other hand, forgetting this ordering and
the distinction between ``inputs'' and ``outputs'' gives a natural map
$\R(p,q)\to\R(p+q)$; in particular, if $\Gamma\in\R(n)$, then two pre-images of
$\Gamma$ in $\R(n,0)$ may only differ by a permutation of the indices on
the inputs. Thus, we can regard all the legs of $\Gamma$ as inputs, and
define a linear map $Z_A(\Gamma)$ on the sub-space of $\Perm{n}$-invariant
vectors of $V\tp{n}$:
\begin{equation*}
  Z_{A}(\Gamma)\colon \bigl(V\tp{n}\bigr)^{\Perm{n}}\to \setC[x_*],
  \qquad
  \Gamma\in\R(n).
\end{equation*}
Note, in particular, that the amplitude of a closed ribbon graph is a well
defined polynomial: indeed, the canonical map $\R(0,0)\to\R(0)$ is an
identification.

\begin{definition}\label{dfn:expectation-value}
  Let \(\Psi\) be any ribbon graph. Its \emph{expectation
    value} is the formal series in the variables \(x_*\):
  \begin{equation*}
    \label{eq:expectation}
    \langle\!\langle \Psi\rangle\!\rangle_{A}^{} 
    := \sum_{\Gamma \in \R_{\Psi}^{}(0)}
    \frac{Z_{A}^{}(\Gamma)}{\card{\Aut{\Gamma}}}\,.
  \end{equation*}
\end{definition}
If the graph \(\Psi\) has \(n\) legs, the function 
\begin{math}
  v\mapsto Z_{A}^{}(\Psi)\bigl(v^{\otimes n}\bigr)
\end{math}
is a well defined polynomial map $V\to \setC[x_*]$. Therefore, it is
integrable on $V_\setR$ with respect to the normalized Gaussian measure
\begin{equation*}
  \ud\mu := \frac{\displaystyle{\exp \left\{-\frac{1}{2}b(v,v) \right\}\ud
      v}}{\displaystyle{\int_{V_\setR}\exp \left\{-\frac{1}{2}b(v,v) \right\}
      \ud v}}.
\end{equation*}
\begin{definition}\label{definition:potentialalgebra}
  The formal series 
  \begin{equation*}
    S_A(x_*):=\sum_{k=1}^\infty x_k
    \frac{T_k(v\tp{k})}{k}
  \end{equation*}
  is called the \emph{potential} of the cyclic algebra $A$.
\end{definition}

With the above notations, we have the following fundamental formula.
\begin{proposition}[Feynman-Reshetikhin-Turaev]
  \label{prop:FRT}
  For any ribbon\break graph $\Psi\in\R(n)$, the following asymptotic
  expansion holds:
  \begin{equation}
    \label{eq:FRT1}
    \expval{\Psi}_{A}^{} = \int_{V_\setR} \frac{Z_{A}
      (\Psi)(v\tp{n})}{\card{\Aut\Psi}}
    \exp S_A(x_*) \ud\mu(v) 
  \end{equation}
\end{proposition}
For a proof, see \cite[Theorem 3.6. and Formula
3.3]{fiorenza-murri;feynman} and \cite{fiorenza;integration-groupoids}.

In particular, we have:
\begin{align*}
  \langle\!\langle \emptyset\rangle\!\rangle_{A}^{}&=
  \int_{V_\setR}
  \exp S_A(x_*)
  \ud\mu(v)
  \\
  \langle\!\langle{\sf v}_n\rangle\!\rangle_{A}^{}&=
  \int_{V_\setR}\frac{T_n(v\tp{n})}
  {n}
  \exp
  S_A(x_*)\ud\mu(v)
  \\
  \langle\!\langle \coprod_{i=1}^n {{\sf v}_i}^{\coprod
    m_i}
  \rangle
  \!\rangle_{A}^{}
  &=
  \int_{V_\setR}\frac{\prod_{i=1}^n{(T_i(v\tp{i}))}^{m_i}}
  {\prod_{i=1}^n i^{m_i}m_i!}
  \exp S_A(x_*) 
  \ud\mu(v)
\end{align*}
where $\coprod_{i=1}^n {{\sf v}_i}^{\coprod m_i}$ denotes the disjoint union of
$m_1$ copies of ${\sf v}_1$, $m_2$ copies of ${\sf v}_2$, \dots and
$m_n$ copies of ${\sf v}_n$.
\begin{definition}\label{definition:partfunalgebra}
  The formal series
  \begin{equation*}
    Z_A(x_*):=\expval{\emptyset}_{A}^{},
    \qquad 
    F_A(x_*):=\log Z_A(x_*),
  \end{equation*}
  are called, respectively, the \emph{partition function} and the
  \emph{free energy} of the cyclic algebra $A$.
\end{definition}
\begin{remark}\label{rem:alpha}
By definition, the partition function is the weighted sum of the
amplitudes of all closed ribbon graphs with only ordinary vertices; a
standard combinatorial argument (see
\cite{bessis-itzykson-zuber;graphical-enumeration}) proves that the free
energy can be
written as the weighted sum of the amplitudes of all \emph{connected}
closed ribbon graphs with only ordinary vertices.
\end{remark}
\begin{remark}\label{rem:beta}
Note that 
\begin{equation}
  \label{eq:d1Z}
  \frac{\del}{\del x_n} \expval{\emptyset}_{A}^{} 
  =\expval{{\sf v}_n}_{A},
\end{equation}
and, more in general,
\begin{equation}
  \label{eq:dnZ}
  \frac{\del^{m_1+\dots +m_n} }{{\del x_1}^{m_1} 
    \cdots {\del x_n}^{m_n}} \expval{\emptyset}_{A}^{}
  = {m_1!\cdots m_n!} \cdot \expval{\coprod_{i=1}^n {{\sf
      v}_i}^{\coprod m_i}}_{A},
\end{equation}
that is, derivatives of $\expval{\emptyset}_{A}^{}$ can be written as
expectation values of (disjoint unions of) special vertices.
\end{remark}
\subsection{Ribbon graphs with colored edge-sides}
\label{sec:coledgesides}
Let \((V,b,T_1,T_2,\dots)\) be a cyclic algebra and assume that \(V\)
has a decomposition
\begin{math}
  V = \bigoplus_{\xi,\eta\in I} V_{\xi,\eta}
\end{math}
where \(V_{\xi,\eta}\) and \(V_{\eta,\xi}\) are dual subspaces with
respect to the pairing $b$. Since any edge of a ribbon graph has two
(distinct) sides, it is meaningful to consider ribbon graphs with edge
sides colored by elements of $I$. We introduce the following new
graphical calculus element on~$V$:
\begin{equation*}
  Z_{A,I}\left({\xy
      (0,-1);(0,1)**\dir{-},(0,-1.2)*{1_{\text{\rm
            in}}},(0,1.2)*{1_{\text{\rm
            out}}},(-.2,0)*{\xi},(0.2,0)*{\eta}\endxy}\right):=
  \pi_{\xi,\eta}\colon
  V\to V\,,   
\end{equation*}
with \(\pi_{\xi,\eta}\) the orthogonal projection on the
subspace \(V_{\xi,\eta}\).
Since
\begin{align*}
  \xy*!LC\xybox{
    \vloop~{(0,0.5)}{(1,0.5)}{(0,-0.5)}{(1,-0.5)},(0,-0.7)*{
      x} ,(1,-0.7)*{y}},(-0.1,-.2)*{\xi},(0.3,-.2)*{\eta}
  \endxy=b({\pi_{\xi,\eta}(x)},{y})=
  b({x},{\pi_{\eta,\xi}(y)})
  =\xy*!LC\xybox{
    \vloop~{(0,0.5)}{(1,0.5)}{(0,-0.5)}{(1,-0.5)},(0,-0.7)*{
      x} ,(1,-0.7)*{y}},(.9,-.2)*{\eta},(1.4,-.2)*{\xi}
  \endxy,
\end{align*}
then graphical calculus extended to ribbon graphs with colored edge
sides is well defined; that is, a graphical calculus functor $Z_{A,I}$
is defined for any cyclic algebra $A$ and a decomposition into
subspaces indexed by $I$ as above; moreover, $Z_{A,I}$ enjoys the
properties listed in \prettyref{prop:gc1}.

From
\begin{math}
\id_V=\bigoplus_{\xi,\eta\in I}
\pi_{\xi,\eta}
\end{math}
we obtain the graphical identity
\begin{equation*}
Z_{A}\left({\xy
      (0,-1);(0,1)**\dir{-},(0,-1.2)*{1_{\text{\rm
            in}}},(0,1.2)*{1_{\text{\rm
            out}}}\endxy}\right)=\bigoplus_{\xi,\eta\in I}
Z_{A,I}\left({\xy
      (0,-1);(0,1)**\dir{-},(0,-1.2)*{1_{\text{\rm
            in}}},(0,1.2)*{1_{\text{\rm
            out}}},(-.2,0)*{\xi},(0.2,0)*{\eta}\endxy}\right),
\end{equation*}
so the amplitude of a ribbon graph $\Gamma$ is expanded into the sum of
amplitudes of ribbon graphs obtained by coloring the edge sides of
$\Gamma$ with colors in the set $I$, in all possible ways.

\section{The 't~Hooft-Kontsevich model}
\label{sec:tHooft}

The space $M_N(\setC)$ of $N\times N$ complex matrices  has a natural
Hermitian inner product
\begin{equation*} 
  \inner{X}{Y}:=\tr(X^*Y),
\end{equation*}
which induces the standard Euclidean inner product \(\inner{X}{Y} =
\tr(XY)\) on the real subspace of Hermitian matrices
\begin{equation*}
  \Hermitian[N] :=\{X\in M_N(\setC) | X^*=X\}.
\end{equation*}

For any positive definite $N\times N$ Hermitian matrix \(\Lambda\), we
can define a new Euclidean inner product on \(\Hermitian[N]\) as
\begin{equation*}
  \inner[\Lambda]{X}{Y} := 
  {\textstyle\onehalf} \bigl( \tr(X\Lambda Y) + \tr(Y\Lambda X) \bigr).
\end{equation*}

The complexification ${\Hermitian[N]}\otimes \setC$ is canonically
isomorphic to the space $M_N(\setC)$ of $N \times N$ complex matrices,
so the pairing $\inner[\Lambda]{-}{-}$ induces a non-degenerate
symmetric bilinear form on $M_N(\setC)$. Define cyclic tensors
$T_k:M_N(\setC)\tp{k} \to \setC$ by
\begin{equation*}
  T_{k}(X_1\otimes X_2\otimes\cdots\otimes X_{k}) := 
  \tr(X_1\cdots X_k)
\end{equation*} 
The tensors \(T_k\) together with the pairing $b :=
\inner[\Lambda]{-}{-}$ give a cyclic algebra structure on the space of
$N\times N$ complex matrices; denote $Z_\Lambda$ the graphical calculus
functor induced by this cyclic algebra structure on the category of
ribbon graphs (see \prettyref{prop:gc1}).

We ought to compute $Z_{\Lambda}(\Gamma)$ for all the generators,
i.e., for \(k\)-valent vertices with incoming legs and for the bent
edge with outgoing legs. Denote by $\{E_{ij}\}$ the canonical basis of
$M_N(\setC)$; it is immediate to reckon:
\begin{equation}
  \label{eq:TR1}
Z_\Lambda\left(\xy*!LC\xybox{
\vloop~{(0,-0.5)}{(1,-0.5)}{(0,0.5)}{(1,0.5)},(0,0.7)*{\,\,
1^{\textrm{out}}} ,(1,0.7)*{\,\, 2^{\textrm{out}}}}
\endxy\right)
  = \sum_{i,j} \frac {2} {\Lambda_i + \Lambda_j}
  E_{ij} \otimes E_{ji},
\end{equation}
\begin{multline}
  \label{eq:TR2}
  Z_{\Lambda}\left(
    \xyc\rgvertex5\loose1\loose2\loose3\loose4\missing5%
,(0,1.6)*{\,\,1^{\textrm{in}}},(-1.3,.5)*{\,\,2^{\textrm{in}}}%
,(1.5,0.5)*{\,\,k^{\textrm{in}}},(-.8,-1.2)*{\,\,3^{\textrm{in}}}
\endxyc
  \right) (E_{j_k i_1} \otimes E_{j_1 i_2
    } \otimes \cdots \otimes
  E_{ j_{k-1} i_k})
  \\
  = x_kT_k (E_{j_k i_1} \otimes E_{j_1 i_2} \otimes \cdots \otimes
  E_{j_{k-1} i_k})
  \\
  =x_k \cdot \delta_{i_1j_1}
  \delta_{i_2j_2}\cdots \delta_{i_k j_k},
\end{multline}
and
\begin{multline}
  \label{eq:TR3}
  Z_{\Lambda}\left(
    \xyc\rgvertex[\!]{5}\loose1\loose2\loose3\loose4\missing5%
,(0,1.6)*{\,\,1^{\textrm{in}}},(-1.3,.5)*{\,\,2^{\textrm{in}}}%
,(1.5,0.5)*{\,\,k^{\textrm{in}}},(-.8,-1.2)*{\,\,3^{\textrm{in}}}
\endxyc
  \right) (E_{j_k i_1} \otimes E_{j_1 i_2
    } \otimes \cdots \otimes
  E_{ j_{k-1} i_k})
  \\
  =\delta_{i_1j_1}
  \delta_{i_2j_2}\cdots \delta_{i_k j_k}.
\end{multline}

For any \(i,j\in\{1,\dots,N\}\), let
\begin{math}
  M_{ij}=\setC\cdot E_{ij}\,;
\end{math}
then 
\begin{math}
M_N(\setC)=\bigoplus_{i,j}M_{ij}\,,
\end{math}
and \(M_{ij}\) is the dual of \(M_{ji}\) with respect to the pairing
\({\inner{-}{-}}_\Lambda\), so $Z_\Lambda$ is actually a graphical
calculus for ribbon graphs with edge sides colored with indices from
$I = \{1,\dots,N\}$.

Decorating the sides of a leg with the indices $i,j$ from
$\{1,\dots,N\}$ is coherent with the convention of writing $i,j$ near
an endpoint to denote evaluation at the basis element $E_{ij}$.
Formulas \eqref{eq:TR1}--\eqref{eq:TR3} can therefore be rewritten as:
\begin{equation*}
  {\xyc\vloop-
    ,(0.05,0.05)*\txt{\({}_i\ {}_j\)},(1,0.05)*\txt{\({}_j\ {}_i\)}
    \endxyc}
  = \frac {2} {\Lambda_i + \Lambda_j} E_{ij}\otimes E_{ji},
\end{equation*}
\begin{equation*}
  \label{eq:EV1}
  {\xyc%
    \rgvertex5%
    \fence3{j_1}{i_2}%
    \fence4{j_2}{\ldots}%
    \missing5%
    \fence1{\ldots}{i_{k}}%
    \fence2{j_k}{i_1}%
    \endxyc}
= \begin{cases}
  x_k, &\text{if $i_l = j_l$, for all $l$,}
  \\
  0,   &\text{otherwise,}
\end{cases}
\end{equation*}
and
\begin{equation*}
  \label{eq:EV2}
  {\xyc%
    \rgvertex[\!]5%
    \fence3{j_1}{i_2}%
    \fence4{j_2}{\ldots}%
    \missing5%
    \fence1{\ldots}{i_{k}}%
    \fence2{j_k}{i_1}%
    \endxyc} 
= \begin{cases}
  1, &\text{if $i_l = j_l$, for all $l$,}
  \\
  0,   &\text{otherwise.}
\end{cases}
\end{equation*}
That is, according to \prettyref{eq:TR2}, a vertex gives a non-zero
contribution if and only if sides belonging to the same hole are
decorated with the same index.
 
Summing up, for any closed ribbon graph $\Gamma$, one has:
\begin{equation}\label{eq:like-KMI-rhs}
  Z_{\Lambda}(\Gamma)=\prod_{k=1}^\infty{x_k}^{m_k}\sum_{c}\prod_
  {l\in\Edges{\Gamma}}\frac{2}{\Lambda_{c(l^+)}+\Lambda_{c(l^-)}},
\end{equation}
where $m_k$ is the number of ordinary $k$-valent vertices of $\Gamma$, $c$
ranges in the set of all maps $\Holes{\Gamma}\to\{1,\dots,N\}$, and $l^\pm$ are
the two (not necessarily distinct) holes $l$ belongs to.

\subsection{The 't~Hooft-Kontsevich model}
The right hand side of equation \prettyref{eq:like-KMI-rhs} is similar
to the right hand side of the Kontsevich's Main Identity
\prettyref{eq:KMI2}; indeed, we can tie graphical calculus to
Kontsevich' results as follows.
\begin{definition}
  Denote $Z_{\Lambda,s_*}$ the functor obtained from the graphical calculus
  $Z_\Lambda$ by taking:
  \begin{equation}\label{eq:xs}
    x_k=\begin{cases}
      0,                                        &\text{if $k$ is even,}
    \\
      -\sqrt{-1}\left(-\frac{1}{2}\right)^rs_r, &\text{if $k=2r+1$ is odd}.
    \end{cases}
  \end{equation}
  The resulting graphical calculus $Z_{\Lambda,s_*}$ is called the
  't~Hooft-Kontsevich model.
\end{definition}
The partition function (see \prettyref{definition:partfunalgebra}) of
the  't~Hooft-Kontsevich model is the formal series in
the variables $s_*$:
\begin{equation}\label{eq:partHK}
\langle\!\langle\emptyset\rangle\!\rangle_{\Lambda,s_*}^{}:=
\sum_{\Gamma\in\R_\emptyset(0)}\frac{Z_{\Lambda,s_*}(\Gamma)}{\card{
\Aut\Gamma}}.
\end{equation}
It is an asymptotic expansion of the matrix integral
\begin{equation*}
   \int_{\Hermitian[N]} \exp S(s_*;X)
  \ud\mu_\Lambda^{}(X)\,,
\end{equation*}
where $S(s_*,X)$ is the potential of the 't~Hooft-Kontsevich model,
given by
\begin{equation*}
  S(s_*,X) = -\sqrt{-1} \sum_{j=0}^\infty (-1/2)^j s_j
  \frac{\tr X^{2j+1}} {2j+1},
\end{equation*}
and $\ud\mu_\Lambda^{}$ is the Gaussian measure on $\Hermitian[N]$
induced by the pairing $\inner[\Lambda]{-}{-}$ --- see
\prettyref{definition:potentialalgebra} and \prettyref{prop:FRT}.

To compute the partition function
$\langle\!\langle\emptyset\rangle\!\rangle_{\Lambda,s_*}^{}$ we have
to compute the amplitudes of closed ribbon graphs with only
ordinary vertices. 
\begin{definition}
We shall say that a closed ribbon graph has
combinatorial type $m_*$ if, for any $i$, it has exactly $m_i$
ordinary vertices of valence $2i+1$, and no special vertices.
\end{definition}
\begin{lemma}
  In the $N$-dimensional 't~Hooft-Kontsevich model, for any closed
  ribbon graph $\Gamma$ of combinatorial type $m_*$ with $n$ holes, the
  following formula holds:
  \begin{equation}\label{eq:lookslike}
    Z_{\Lambda,s_*}(\Gamma)=(-1)^n\prod_{r=0}^\infty\left(
      \frac{s_r}{2^r}\right)^{m_r}
    \sum_{c}\prod_{l\in\Edges{\Gamma}}\frac{2}{\Lambda_{c(l^+)}+
\Lambda_{c(l^-)}},
  \end{equation}
where $c$ ranges in the set of all maps $\Holes{\Gamma}\to\{1,\dots,N\}$,
and $l^\pm$ are
the two (not necessarily distinct) holes $l$ belongs to.
\end{lemma}
\begin{proof}
  By formula \prettyref{eq:like-KMI-rhs} we immediately obtain:
  \begin{equation*}
    Z_{\Lambda,s_*}(\Gamma)
    = (-\sqrt{-1})^{\card{\Vertices{\Gamma}}}(-1)^{\sum_{j=0}^\infty
      jm_j}\prod_{r=0}^\infty
    \left(
      \frac{s_r}{2^r}\right)^{m_r}
    \sum_{c}\prod_
    {l\in\Edges{\Gamma}} 
    \frac{2}{\Lambda_{c(l^+)}+\Lambda_{c(l^-)}}.
  \end{equation*}
  Thus, we just need to prove $(-\sqrt{-1})^{\card{\Vertices{\Gamma}}}
  \cdot (-1)^{\sum_{j=0}^\infty jm_j} = (-1)^n$; the ribbon graph
  $\Gamma$ satisfies the combinatorial relations:
  \begin{align*}
    \card{\Vertices{\Gamma}} &= {\textstyle\sum\nolimits_j} m_j,
    \\
    2\card{\Edges{\Gamma}} &= {\textstyle\sum\nolimits_j} (2j+1)m_j=
    2({\textstyle\sum\nolimits_j} jm_j)+\card{\Vertices{\Gamma}},
    \\
    \card{\Vertices{\Gamma}} - \card{\Edges{\Gamma}} + n &= 
    \card{\Vertices{\Gamma}} - \card{\Edges{\Gamma}} +
    \card{\Holes{\Gamma}} = \chi\bigl({S(\Gamma)}\bigr) \equiv 0 \pmod 2,
  \end{align*}
  where $S(\Gamma)$ is the Riemann surface associated to $\Gamma$ and $\chi$
  denotes the Euler-Poincar{\'e} characteristic.  Therefore,
  \begin{equation*}
    (-1)^n = \sqrt{-1}^{2n} = \sqrt{-1}^{-\card{\Vertices{\Gamma}}}
    \sqrt{-1}^{2\sum_j jm_j}
    = (-\sqrt{-1})^{\card{\Vertices{\Gamma}}} (-1)^{\sum_j jm_j}.
  \end{equation*}
\end{proof}

\begin{proposition}
  \label{prop:Z=Z}
  The partition function $Z(s_*;t_*)$ of combinatorial intersection
  numbers and the partition function
  $\langle\!\langle\emptyset\rangle\!\rangle_{\Lambda,s_*}^{}$ of the
  't~Hooft-Kontsevich model are related by:
  \begin{equation*}
    Z(s_*;t_*)\vert_{t_*(\Lambda)} =
\langle\!\langle\emptyset\rangle\!\rangle_{\Lambda,s_*}^{}.
  \end{equation*}
\end{proposition}
\begin{proof}
  The statement is clearly equivalent to proving that the same relation
  holds between the free energies, i.e.,
  \begin{equation*}
    F(s_*;t_*)\vert_{t_*(\Lambda)} =
\log\langle\!\langle\emptyset\rangle\!\rangle_{\Lambda,s_*}^{}.
  \end{equation*}
  By \prettyref{rem:alpha} and formula \prettyref{eq:lookslike} we
immediately get:
  \begin{equation*}
\log\langle\!\langle\emptyset\rangle\!\rangle_{\Lambda,s_*}^{}=\sum_{m_*,n}
\sum_{\Gamma,c}
    \frac{(-1)^n} {\card{\Aut{\Gamma}}}
    \prod_{k=0}^\infty \left(\frac{s_k}{2^k}\right)^{m_k} 
    \prod_{l\in\Edges{\Gamma}} \frac{2} 
    {\Lambda_{c(l^+)} + \Lambda_{c(l^-)}},
  \end{equation*}
  where $\Gamma$ ranges over \emph{connected} closed numbered ribbon graphs
  of combinatorial type $m_*$ with \emph{only ordinary vertices}, and
  $c$ is a coloring of $\Holes{\Gamma}$ with colors $\{1,\dots,N\}$.  Now,
  any $c: \Holes{\Gamma} \to \{1, \ldots, N\}$ factors in $n!$ ways as $j \circ h$ where
  $h$ is a bijection $h: \Holes{\Gamma} \to \{1, \ldots, n\}$ and $j$ is a map $j:
  \{1, \ldots, n\} \to \{1, \ldots, N\}$, so we can rewrite the above equation as:
  \begin{equation*}
    \log\langle\!\langle\emptyset\rangle\!\rangle_{\Lambda,s_*}^{}=
\sum_{m_*,n} \sum_{\Gamma, h, j}
    \frac{(-1)^n}{n!\card{\Aut{\Gamma}}}
    \prod_{k=0}^\infty \left(\frac{s_k}{2^k}\right)^{m_k}
    \prod_{l\in\Edges{\Gamma}} \frac{2} 
    {\Lambda_{(j\circ h)(l^+)} + \Lambda_{(j\circ h)(l^-)}}.
  \end{equation*}

  The group $\Aut\Gamma$ acts on the sets $\Edges{\Gamma}$ and  
  $\Holes{\Gamma}$; in particular, the second action induces an action of
  $\Aut\Gamma$ on the set 
  \begin{equation*}
    \text{Num}(\Gamma) := \bigl\{ h : \Holes{\Gamma} \xrightarrow{\sim} 
    \{1, \ldots, n\} \bigr\}.
  \end{equation*}
  It is immediate to check that, if $h_1$ and $h_2$ are in the same orbit
  with respect to the action of $\Aut\Gamma$, then
  \begin{equation*}
    \prod_{l\in\Edges{\Gamma}}\frac{2}{\Lambda_{(j\circ
        h_1)(l^+)}+
      \Lambda_{(j\circ h_1)(l^-)}}
    =
    \prod_{l\in\Edges{\Gamma}}\frac{2}{\Lambda_{(j\circ
        h_2)(l^+)}+
      \Lambda_{(j\circ h_2)(l^-)}},
  \end{equation*}
  so that:
  \begin{multline*}
\log\langle\!\langle\emptyset\rangle\!\rangle_{\Lambda,s_*}^{}=\sum_{m_*,n}
\sum_{\Gamma,[h], j} \Biggl\{
    \frac{(-1)^n \card{\text{orbit of $h$}}} {n!\card{\Aut{\Gamma}}}
    \prod_{k=0}^\infty \left(\frac{s_k}{2^k}\right)^{m_k} \times
    \\
    \prod_{l\in\Edges{\Gamma}} \frac{2} 
    {\Lambda_{(j\circ h)(l^+)} + \Lambda_{(j\circ h)(l^-)}}
    \Biggr\},
  \end{multline*}  
  and this time we take one representative $h$ from each orbit $[h]$
  in $\text{Num}(\Gamma)$.
  
  The isotropy subgroup of any $h \in \text{Num}(\Gamma)$ is $\Aut(\Gamma, h)$;
  therefore, the cardinality of the orbit of $h$ is $\card{\Aut\Gamma} /
  \card{\Aut (\Gamma, h)}$. Moreover, the numbered ribbon graphs $(\Gamma,h_1)$
  and $(\Gamma,h_2)$ are isomorphic if and only if $h_1$ and $h_2$ are in
  the same orbit with respect to the action of $\Aut\Gamma$ on
  $\text{Num}(\Gamma)$, therefore:
  \begin{multline*}
    \log\langle\!\langle\emptyset\rangle\!\rangle_{\Lambda,s_*}^{}=\sum_{m_*,n}
    \sum_{(\Gamma,h), j} \Biggl\{
    \frac{(-1)^n} {n!\card{\Aut{(\Gamma,h)}}}
    \prod_{k=0}^\infty \left(\frac{s_k}{2^k}\right)^{m_k} \times
    \\
    \prod_{l\in\Edges{\Gamma}} \frac{2}
    {\Lambda_{(j\circ h)(l^+)} + \Lambda_{(j\circ h)(l^-)}}
    \Biggr\},
  \end{multline*}
  where $(\Gamma,h)$ ranges over the set of isomorphism classes of
  \emph{connected closed numbered} ribbon graphs of combinatorial type
  $m_*$ with $n$ holes.

  Finally, by Kontsevich' Main Identity \eqref{eq:KMI2} we get:
  \begin{align*}
    \log\langle\!\langle\emptyset\rangle\!\rangle_{\Lambda,s_*}^{} &=
    \sum_{n,m_*,\nu_*} \frac{1}{n!}  s_*^{m_*}\langle
    \tau_{\nu_1}\tau_{\nu_2}\cdots\tau_{\nu_n}\rangle_{m_*,n} \sum_{j}
    \prod_{i=1}^n \frac{-(2\nu_i-1)!!}{{\Lambda_{j(i)}}^{2\nu_i+1}}
    \intertext{where $j\colon\{1, \ldots, n\} \to \{1, \ldots, N\}$,}
    \\
    &=\sum_{n;m_*;\nu_*}\frac{1}{n!}  s_*^{m_*}\langle
    \tau_{\nu_1}\tau_{\nu_2}\cdots\tau_{\nu_n}\rangle_{m_*,n}
    \prod_{i=1}^n\bigl(-(2\nu_i+1)!!\tr\Lambda^{-(2\nu_i+1)}\bigr)
    \\
    &=\sum_{n;m_*;\nu_*}\frac{1}{n!}  s_*^{m_*}\langle
    \tau_{\nu_1}\tau_{\nu_2}\cdots\tau_{\nu_n}\rangle_{m_*,n}
    t_{\nu_1}(\Lambda)\cdots t_{\nu_n}(\Lambda)
    \\
    &=F(s_*;t_*)\bigr\rvert_{t_*(\Lambda)}.
  \end{align*}
\end{proof}

\section{Witten's formula for derivatives}
\label{sec:derivatives}

It has been remarked by Witten \cite{witten;kontsevich-model} that
first order derivatives of the partition function $Z(t_*)$ are related
to expectation values in the $(N+1)$-dimensional
't~Hooft-Kontsevich model of ribbon graph with one ``distinguished'' hole.
We shall use graphical calculus for ribbon graphs with sides of two
colors in order to present a version of Witten's argument suitable for
application to the partition function $Z(s_*; t_*)$.

\subsection{The $(N+1)$-dimensional 't~Hooft-Kontsevich model}
\label{sec:N+1}

Let \(z\) be a \emph{real} positive variable, and consider the
$(N+1)$-dimensional 't~Hooft-Kontsevich model $Z_{z\oplus\Lambda,
  s_*}$ defined by the diagonal
matrix 
\begin{equation*}
  z\oplus \Lambda = 
  \begin{pmatrix}
    z       & 0\\
    0 & \Lambda 
  \end{pmatrix}.
\end{equation*}
Let \(\{E_{ij}\}_{i,j=0,\dots,N}\) be the canonical basis for
\(M_{N+1}(\setC)\), and define
\begin{align*}
&M^{N+1}_{\Lambda,\Lambda}:=\text{span}(E_{ij})_{i,j>0} \sim M_N(\setC),
&& M^{N+1}_{z,\Lambda}:=\text{span}(E_{0j})_{j>0},
\\
&M^{N+1}_{\Lambda,z}:=\text{span}(E_{i0})_{i>0},
&& M^{N+1}_{z,z}:=\setC\cdot E_{00}\,.
\end{align*}
Then we have the following decomposition:
\begin{equation*}
M_{N+1}(\setC)=\bigoplus_{\xi,\eta\in\{\Lambda,z\}} M^{N+1}_{\xi,\eta}\,,
\end{equation*}
therefore a graphical calculus for ribbon graphs with sides colored
with the two colors \(\Lambda\) and \(z\) is defined on
\(M_{N+1}(\setC)\), extending the Feynman rules for the
`t~Hooft-Kontsevich model.

Note that a ribbon graph with only \(\Lambda\)-decorated edges is naturally
identified to a graphical element for the \(N\)-dimensional
`t~Hooft-Kontsevich model; for this reason, we will omit $\Lambda$ from the
decoration of edges in the displayed diagrams of this paper.

Moreover, we will put a \(z\) in the middle of an hole to mean that
all the edge-sides of that hole are \(z\)-decorated, e.g.,
\begin{equation*}
{\xy
\rghole{3}\loose{1}\loose{2}\loose{3},(0,0)*{z}%
\endxy}\qquad
{:=}
\qquad{\xy
\rghole{3}\loose{1}\loose{2}\loose{3}%
,(-0.33,0.13)*{z},
,(0.33,0.13)*{z},
,(0,-0.35)*{z},
\endxy}
\end{equation*}

By the above definitions it is immediate to compute propagators in
the \((N+1)\)-dimensional `t~Hooft-Kontsevich model:
\begin{align}
  \label{eq:casimirs-n-piu-uno}
  Z_{z\oplus\Lambda, s_*}\left(\xy*!LC\xybox{
      \vloop~{(0,-0.5)}{(1,-0.5)}{(0,0.5)}{(1,0.5)},(0,0.7)*{\,\,
        1^{\textrm{out}}} ,(1,0.7)*{\,\, 2^{\textrm{out}}}}%
    ,(.93,-.6)*{z},(.93,-.93)*{z} \endxy\right)&=\frac{1}{z} E_{00}
  \otimes E_{00},
  \\
  Z_{z\oplus\Lambda, s_*}\left(\xy*!LC\xybox{
      \vloop~{(0,-0.5)}{(1,-0.5)}{(0,0.5)}{(1,0.5)},(0,0.7)*{\,\,
        1^{\textrm{out}}} ,(1,0.7)*{\,\, 2^{\textrm{out}}}}%
    ,(.93,-.6)*{z} \endxy\right)&=\sum_{i=1}^N \frac{2}{z + \Lambda_i}
  E_{i0} \otimes E_{0i},
  \\
  Z_{z\oplus\Lambda, s_*}\left(\xy*!LC\xybox{
      \vloop~{(0,-0.5)}{(1,-0.5)}{(0,0.5)}{(1,0.5)},(0,0.7)*{\,\,
        1^{\textrm{out}}} ,(1,0.7)*{\,\, 2^{\textrm{out}}}}%
    ,(.93,-.93)*{z}
    \endxy\right)&= \sum_{i=1}^N \frac{2}{\Lambda_i + z} E_{0i}
  \otimes E_{i0},
  \\
  Z_{z\oplus\Lambda,
    s_*}\left(\xy*!LC\xybox{
      \vloop~{(0,-0.5)}{(1,-0.5)}{(0,0.5)}{(1,0.5)},(0,0.7)*{\,\,
        1^{\textrm{out}}} ,(1,0.7)*{\,\, 2^{\textrm{out}}}}
    \endxy\right)&= \sum_{i,j=1}^N \frac{2}{\Lambda_i
    + \Lambda_j} E_{ij}
  \otimes E_{ji}.
\end{align}

Since amplitudes of vertices are null if two consecutive sides are not
decorated with the same color, then, reasoning as in
\prettyref{sec:tHooft}, one finds that the amplitude
\(Z_{z\oplus\Lambda,s_*}^{}(\Gamma)\) of a ribbon graph $\Gamma$ in the
$(N+1)$-dimensional 't~Hooft-Kontsevich model can be diagrammatically
written as a sum of copies of \(\Gamma\) with some of the holes
decorated by the variable \(z\).

\subsection{Witten's formula}
The machinery is now in place to prove Witten's formula.
\begin{proposition}
For any $k\geq0$, the following identity holds:
\begin{equation}\label{eq:dZst}
  \left.\frac {\del Z(s_*;t_*)} {\del t_k}\right\rvert_{t_*(\Lambda)}
  = -\frac{1}{(2k-1)!!}  
  \Coeff_{z}^{-(2k + 1)}\! \Biggl(
    \sum_{\Gamma\in\R_\emptyset^{[1]}(0)}
    \frac{Z_{z\oplus\Lambda,s_*}(\Gamma)} {\card{\Aut{\Gamma}}}
    \Biggr),  
\end{equation}
where $\R_\emptyset^{[1]}(0)$ denotes the set of isomorphism classes
of closed ribbon graphs with only ordinary vertices and exactly
one $z$-decorated hole.
\end{proposition}
\begin{proof}
  Since a $z$-decorated hole lies in one connected component of the
  ribbon graph, by the usual combinatorial argument we can reduce to
connected ribbon
  graphs, i.e., the statement is equivalent to:
\begin{equation}\label{eq:dFst}
\left.\frac {\del F(s_*;t_*)} {\del
t_k}\right\rvert_{t_*(\Lambda)}=
-\frac{1}{(2k-1)!!}  \Coeff_{z}^{-(2k + 1)}\! \left(
\sum_{\Gamma}
\frac{Z_{z\oplus\Lambda,s_*}(\Gamma)}{\card{\Aut{\Gamma}}}
    \right),    
\end{equation}
with $\Gamma$ ranging in the set of isomorphism classes of connected
closed ribbon graphs with exactly one $z$-decorated hole and with only
generic vertices.  To prove equation \prettyref{eq:dFst}, note that
the Kontsevich' Main Identity (\prettyref{prop:KMI2}) for graphs with
\(n+1\) holes numbered from \(0\) to \(n\) gives:
  \begin{multline*}
    \sum_{k} \left(s_*^{m_*} \sum_{\nu_1, \dots, \nu_n} \langle
      \tau_{k} \tau_{\nu_1} \cdots \tau_{\nu_n} \rangle_{m_*,n+1}
      \prod_{i=1}^n \frac{{(2 \nu_i-1)!!}}{ {\lambda_{i}^{2 \nu_i+1}}}
    \right) \frac{{(2 k-1)!!}}{ {\lambda_0^{2 k+1}}}
    \\
    =\sum_{ (\Gamma,h)} \frac{1}{\card{\Aut{(\Gamma,h)}}}
    \prod_{r=0}^\infty\left( \frac{s_r}{2^r}\right)^{m_r}
    \prod_{l\in\Edges{\Gamma}} \frac{2}{\lambda_{h(l^+)} +
      \lambda_{h(l^-)}},
  \end{multline*}
where $(\Gamma,h)$ ranges in the set of isomorphism classes of closed
connected ribbon graphs of combinatorial type $m_*$, with $n+1$ holes,
numbered from $0$ to $n$.  Therefore,
\begin{multline*}
  \sum_{m_*,\nu_*} s_*^{m_*}\langle \tau_{k}
  \tau_{\nu_1}\cdots\tau_{\nu_n} \rangle_{m_*,n+1}
  \prod_{i=1}^n\frac{(2 \nu_i-1)!!}{\lambda_{i}^{2n_i+1}} =
  \frac{1}{(2k-1)!!} \times
  \\
  \times \Coeff_{\lambda_0}^{-(2k + 1)} \Biggl( \sum_{ (\Gamma,h)}
  \frac{1}{\card{\Aut{(\Gamma,h)}}} \prod_{r=0}^\infty\left(
    \frac{s_r}{2^r}\right)^{m_r}\prod_{l\in\Edges{\Gamma}}
  \frac{2}{\lambda_{h(l^+)}+ \lambda_{h(l^-)}} \Biggr).
\end{multline*}
Now set \(\lambda_i = \Lambda_{j(i)}\), for \(i=0,1, \dots, n\); and sum
over all \(j: \{0,1, \dots, n\} \to \{0,1, \dots,N \}\) such that
$j(0)=0$ and recall that $\Lambda_0=z$ to get:
 \begin{multline*}
   \sum_{m_*,\nu_*,j} s_*^{m_*}\langle \tau_{k}
   \tau_{\nu_1}\cdots\tau_{\nu_n} \rangle_{m_*,n+1}
   \prod_{i=1}^n\frac{(2 \nu_i-1)!!}{\Lambda_{j(i)}^{(2n_i+1)}}
   =\frac{1}{(2k-1)!!}  \times
   \\
   \times 
   \Coeff_{z}^{-(2k + 1)} 
   \Biggl( \sum_{ (\Gamma,h),j} \frac{1}{\card{\Aut{(\Gamma,h)}}}
     \prod_{r=0}^\infty\left( \frac{s_r}{2^r}\right)^{m_r} \!\!
     \prod_{l\in\Edges{\Gamma}} \frac{2}{\Lambda_{j\circ h(l^+)}+
       \Lambda_{j\circ h(l^-)}} \Biggr).
  \end{multline*}
The proof can then be easily concluded by using the fact that
\begin{equation*}
  \frac {\del F(s_*;t_*)} {\del t_k} = \sum_{m_*,\nu_*}
  \frac{1}{n!}\langle 
  \tau_k \tau_{\nu_1} \ldots \tau_{\nu_n} \rangle 
  s_*^{m_*}t_{\nu_1} \cdots t_{\nu_n},
\end{equation*}
and reasoning as in the proof of \prettyref{prop:Z=Z}.
\end{proof}

\subsection{Hole types}
\label{sec:hole-types}
If $\Gamma$ is an element of $\R_\emptyset^{[1]}$, i.e., a closed ribbon
graph with only ordinary vertices and exactly one hole decorated by the
variable $z$, then its $z$-decorated hole can be regarded as a
distinguished sub-diagram of $\Gamma$.
\begin{definition}
A ($z$-decorated) hole type is a ribbon graph with only ordinary
vertices and exactly one $z$-decorated hole, which is minimal with
respect to this property, i.e., such that none of its proper subgraphs
contains the $z$-decorated hole. The set of isomorphism
classes of hole types will be denoted by the symbol ${\mathcal S}$.
\end{definition}
Having introduced this terminology, the previous remark can be
restated as:
\begin{equation*}
\R_\emptyset^{[1]}(0)=\bigcup_{\Gamma\in{\mathcal S}} \R_\Gamma^{[1]}(0),
\end{equation*}
where $\R_\Gamma^{[1]}$ denotes the set of isomorphism classes of
closed ribbon graphs containing the hole type $\Gamma$ as
a distinguished subgraph and having no $z$-decorated hole apart from the
hole of $\Gamma$.  Therefore, we can rewrite \prettyref{eq:dZst} as
\begin{equation}\label{eq:dZst-two}
  \left.\frac {\del Z(s_*;t_*)} {\del t_k}\right\rvert_{t_*(\Lambda)}
  = -\frac{1}{(2k-1)!!}\sum_{\Gamma\in{\mathcal S}}
  \Coeff_{z}^{-(2k + 1)}\! \Biggl(
    \sum_{\Phi\in\R_\Gamma^{[1]}(0)}
    \frac{Z_{z\oplus\Lambda,s_*}(\Phi)} {\card{\Aut{\Phi}}}
    \Biggr).
\end{equation}
By introducing the shorthand notation
\begin{equation}\label{eq:dfn-expval-n+1}
  \expval{\Gamma}_{z\oplus\Lambda,s_*}^{[1]}
  := \sum_{\Phi\in\R^{[1]}_\Gamma(0)}
  \frac{Z_{z\oplus\Lambda}^{}(\Phi)} {\card{\Aut\Phi}}\,,
\end{equation}
where $\Gamma$ is an hole type, Witten's formula for derivatives
finally becomes:
\begin{equation}\label{eq:dZst-twotwo}
  \left.\frac {\del Z(s_*;t_*)} {\del t_k}\right\rvert_{t_*(\Lambda)}
  = -\frac{1}{(2k-1)!!}\sum_{\Gamma\in{\mathcal S}}
  \Coeff_{z}^{-(2k + 1)}\expval{\Gamma}_{z\oplus\Lambda,s_*}^{[1]}
 \end{equation}

\section{Proof of the Main Theorem}
\label{sec:proof}

By the correspondence between Gaussian integrals and expectation values of
graphs (see \prettyref{prop:FRT}), if $\Gamma$ is an hole type with $n$
legs, then
\begin{equation} 
\label{eq:int-1}
\langle\!\langle\Gamma\rangle\!\rangle_{z\oplus\Lambda,s_*}^{[1]} =
\int_{\Hermitian[N]}\frac{Z_{z\oplus\Lambda}^{}(\Gamma)}
{\card{\Aut\Gamma}}\bigl(X^{\otimes n}\bigr) \exp
S(s_*;X)\ud\mu_\Lambda^{}(X)\,.
\end{equation}
Therefore, equation \prettyref{eq:dZst-twotwo} is equivalent to   
\begin{multline}\label{eq:dZst-three}
  \left.\frac {\del Z(s_*;t_*)} {\del t_k}\right\rvert_{t_*(\Lambda)}
  = -\frac{1}{(2k-1)!!} \times \\
  \sum_{\Gamma\in{\mathcal
      S}}\int_{\Hermitian[N]} \negquad \Coeff_{z}^{-(2k +
    1)}\left(\frac{Z_{z\oplus\Lambda,s_*}^{}(\Gamma)}
    {\card{\Aut\Gamma}}\bigl(X^{\otimes n}\bigr)\right)  \times \\
  \exp S(s_*;X)\ud\mu_\Lambda^{}(X);
\end{multline}
so we are interested in the tensors
\begin{math}
  \Coeff_{z}^{-(2k + 1)} {Z_{z\oplus\Lambda},s_*^{}(\Gamma)} /
    {\card{\Aut\Gamma}}(X^{\otimes n}).
\end{math}

\subsection{Special vertices decorated by polynomials}
\label{sec:rpoly}

If $\Gamma$ is a $z$-hole type, then its amplitude $Z_{z\oplus\Lambda;s_*}$ has a
Laurent expansion in powers of $z^{-1}$ as $z\to\infty$, whose coefficients
are polynomials in the \(\Lambda_i\)'s.  Indeed, by the Feynman rules for
the $(N+1)$-dimensional 't~Hooft-Kontsevich model, we have:
\begin{enumerate}
\item each \((2r+1)\)-valent ordinary vertex brings a
factor \(-\U(-1/2)^r\,s_r\);
\item each internal edge bordering the \(z\)-decorated hole on
both sides contributes a factor \(1/z\);
\item the other internal edges  contribute factors
of the form \(2/(z+\Lambda_i)\) for \(i=1,\dots,N\).
\end{enumerate}
In other words, the structure constants of the tensors
$\Coeff_z^{-k}Z_{z\oplus\Lambda;s_*}(\Gamma)$ are polynomials in the
$\Lambda_i$'s. We can therefore graphically represent these tensors
enlarging the class of special vertices by adding special vertices
decorated by polynomials. This is formally done as follows.

Let \(\varphi\) be a polynomial in
\(\setC[\theta_1,\theta_2,\dots,\theta_n]\); we say that the
polynomial \(\varphi\) is cyclically invariant iff it is invariant
with respect to the natural action of the cyclic group \({\mathbb
  Z}/n{\mathbb Z}\) on the coordinates. By the symbol \({\sf
  v}^\varphi_n\) we denote an \(n\)-valent special vertex decorated by
the polynomial \(\varphi\).  We represent graphically these vertices
as:
\begin{align*}
  {\xyc
  \rgvertex[\varphi]{5}%
  \loose{1}%
  \loose{2}%
  \loose{3}%
  \loose{4}%
  \missing{5}
  \endxyc
  }\qquad
 &&&
{\xyc
  \rgvertex[\varphi]{5}%
  \cilia{1}{2}
  \loose{1}%
  \loose{2}%
  \loose{3}%
  \loose{4}%
  \missing{5}
  \endxyc
  \thinspace\thinspace\thinspace
  \text{($n$ edges)}
}
\\
\text{cyclically invariant $\varphi$}
&&&
\text{any polynomial \(\varphi\)}
\end{align*}
The r{\^o}le of the ``$\bigstar$'' mark is precisely to break the cyclical
symmetry of the graphical element. We have to define the Feynman
rules for these new vertices; set
\begin{equation*}
  \xyc
  \rgvertex[\varphi]{5}\genericfences
  \endxyc
  :=
  \varphi(\Lambda_{i_1},\Lambda_{i_2},\dots,\Lambda_{i_n}) \cdot
  \xyc
  \rgvertex[\!]{5}\genericfences
  \endxyc
\end{equation*}
--- this is well defined due to the cyclical invariance of
\(\varphi\) ---, and
\begin{equation*}
  \xyc\rgvertex[\varphi]{5}\genericfences\cilia{2}{3}\endxyc
  :=
  \varphi(\Lambda_{i_1},\Lambda_{i_2},\dots,\Lambda_{i_n})
  \xyc\rgvertex[\!]{5}\genericfences\endxyc,
\end{equation*}
so that the ``$\bigstar$'' tells which indeterminate ---among those
corresponding to indices decorating holes around the vertex--- comes
first. Note that, if $\varphi \in \setC[\theta_1, \ldots, \theta_\nu]$ is non-cyclic and $\nu \neq n$,
then we can nonetheless give ${\sf v}^\varphi_n$ a meaning: indeed, if $\nu <
n$ then the above equation still makes sense; if $\nu > n$ then wrap
around the vertex as many times as needed.  Note that special vertices
decorated by the constant polynomial \(1\) are identified with the
non-decorated special vertices.

Using these notations, the Laurent expansions of
$Z_{z\oplus\Lambda;s_*}(\Gamma)$, for an hole
type $\Gamma$ are easily written as sums over ribbon graphs; we give some
illustrative examples here.

\begin{example}\label{xmp:hole-a}
\begin{align*}
  {\xyc\rghole[z]{3}\fence{1}{i}{j}\fence{2}{j}{k}%
    \fence{3}{k}{i}\endxyc} &= -\U{s_1}^3 \cdot \frac{1}{(z+\Lambda_i) \cdot (z+\Lambda_j)
    \cdot (z+\Lambda_k)}
  \allowbreak \\
  &= -\U {s_1}^3 \cdot 1/z^3 + \U{s_1}^3( {\Lambda_i} +
{\Lambda_j} + {\Lambda_k}) \cdot 1/z^4 + \cdots
  \allowbreak \\
  &= \left( -\U {s_1}^3 \cdot {\xyc\rgvertex[\!]{3}\fence{1}{i}{j}\fence{2}{j}{k}%
      \fence{3}{k}{i}\endxyc} \right) \cdot \frac{1}{z^3} 
  \allowbreak \\
  &\phantom{= -\U {s_1}^3 \cdot 1/z^3}
  + \left( \U {s_1}^3 \cdot
{\xyc\rgvertex[\varphi]{3}\fence{1}{i}{j}\fence{2}{j}{k}%
      \fence{3}{k}{i}\endxyc} \right) \cdot \frac{1}{z^4} 
  + \cdots,
\end{align*}
where
\begin{equation*}
  \varphi(\theta_1, \theta_2, \theta_3) := {\theta_1} +
{\theta_2} +
  {\theta_3}.
\end{equation*}
\end{example}

\begin{example}\label{xmp:hole-b}
  \begin{align*}
    {\xyc
      \rghole[z]{3}%
      \fence3li%
      \fence2hl%
      \fence1jk%
      ,(0,1);(.5,1.4)**\dir{-}?(1)="x",%
      "x"+/r2pt/+/d2pt/*{i};%
      ,(0,1);(-.5,1.4)**\dir{-}?(1)="y",%
      "y"+/l2pt/+/d2pt/*{h};%
      \endxyc} &= \half\U {s_1}^2s_2 \cdot
\frac{1}{(z+\Lambda_i) \cdot (z+\Lambda_l) \cdot
      (z+\Lambda_h)}
    \allowbreak \\
    &= \half\U {s_1}^2s_2 \cdot 1/z^3 - \half\U {s_1}^2s_2 (
{\Lambda_i} +
    {\Lambda_l} + {\Lambda_h} 
    ) \cdot 1/z^4 + \cdots
    \allowbreak \\
    &= \left( \half\U {s_1}^2s_2 \cdot
{\xyc\rgvertex[\!]{5}\fence{1}{i}{j}\fence{2}{j}{k}%
        \fence{3}{k}{h}\fence{4}{h}{l}\fence{5}{l}{i}\endxyc} \right)
    \cdot \frac{1}{z^3} 
    \allowbreak \\
    &\phantom{= \half\U {s_1}^2s_2 }
    - \left( \half\U {s_1}^2s_2 \cdot %
      {\xyc\rgvertex[\varphi]{5}\cilia51%
        \fence{1}{i}{j}\fence{2}{j}{k}%
        \fence{3}{k}{h}\fence4hl\fence5li\endxyc} \right) \cdot
    \frac{1}{z^4} + \cdots
  \end{align*}
where
\begin{equation*}
  \varphi(\theta_1, \theta_2, \theta_3, \theta_4,
  \theta_5) := {\theta_1} + {\theta_4} + {\theta_5}.
\end{equation*}
Note that the coefficient of $z^{-4}$ is not a function of all the
$\Lambda_*$'s around the $z$-decorated hole, so $\varphi$ does not depend on $\theta_2$
and $\theta_3$, and does not exhibit the cyclical invariance found in
\prettyref{xmp:hole-a}. Because of this lack of cyclicity, we use the
``$\bigstar$'' mark.
\end{example}


\begin{example}\label{xmp:hole-d}
\begin{align*}
    &{\xyc*!C\xybox{
        ,(3.1,1.4);(3.05,1.2)**\crv{(3.13,1.5)&(3.5,1.866)%
          &(4,2)&(4.5,1.866)&(5,1)&(4,0)&(3,0.7)}
        ,(3.068,1.0);(2.5,1.4)**\crv{(3.2,1.07)&(3.3,1.09)%
          &(3.3,1.15)&(3.07,1.35)}
        ,(1.9,1.4);(1.95,1.2)**\crv{(1.87,1.5)&(1.5,1.866)%
          &(1,2)&(.5,1.866)&(0,1)&(1,0)&(2,0.7)}
        ,(1.932,1.0);(2.5,1.4)**\crv{(1.8,1.07)&(1.7,1.09)%
          &(1.7,1.15)&(1.93,1.35)}%
        ,(4.75,1);(5.4,1)**\crv{(5.2,1)},%
        ,(4.75,1)*{\bullet},%
        ,(0.25,1);(-0.4,1)**\crv{(-0.2,1)},%
        ,(0.25,1)*{\bullet},%
        ,(3.064,1.0)*{\bullet},(1.938,1.0)*{\bullet},(1.15,1.2)*{z}%
        ,(-0.4,1.2)*{\scriptstyle{i}} ,(5.4,1.2)*{\scriptstyle{j}}
        ,(-0.4,0.8)*{\scriptstyle{i}} ,(5.4,0.8)*{\scriptstyle{j}}
      }\endxyc} = {s_1}^4 \cdot \frac{1}{z(z+\Lambda_i)^2(z+\Lambda_j)^2}
    \allowbreak \\
    &\qquad = {s_1}^4 \cdot 1/z^5 + {s_1}^4 ( -2{\Lambda_i}- 2{\Lambda_j}) \cdot 1/z^6 + \cdots
    \allowbreak \\
    &\qquad = {s_1}^4 \Bigl( {\xyc\rgvertex[\!]6\fence4{i}{i}\endxyc}
    \negqquad\negquad {\xyc\rgvertex[\!]6\fence1{j}{j}\endxyc} \Bigr)
    \cdot \frac{1}{z^5} +
    \allowbreak \\
    &\qquad\qquad + \Biggl[ -2{s_1}^4 \Bigl(
    {\xyc\rgvertex[\varphi]6\fence4{i}{i}\endxyc} \negqquad\!\negqquad
    {\xyc\rgvertex[\!]6\fence1{j}{j}\endxyc} \Bigr) -2{s_1}^4 \Bigl(
    {\xyc\rgvertex[\!]6\fence4{i}{i}\endxyc} \negqquad\!\negqquad
    {\xyc\rgvertex[\varphi]6\fence1{j}{j}\endxyc} \Bigr)\Biggr] \cdot
    \frac{1}{z^6} + \cdots
\end{align*}
where
\begin{math}
  \varphi(\theta_1) = \theta_1.
\end{math}
 \end{example}

\begin{example}\label{xmp:hole-e}
\begin{align*}
  &{\xyc*!C\xybox{
        ,(3.1,1.4);(3.05,1.2)**\crv{(3.13,1.5)&(3.5,1.866)%
          &(4,2)&(4.5,1.866)&(5,1)&(4,0)&(3,0.7)}
        ,(3.068,1.0);(2.5,1.4)**\crv{(3.2,1.07)&(3.3,1.09)%
          &(3.3,1.15)&(3.07,1.35)}
        ,(1.9,1.4);(1.95,1.2)**\crv{(1.87,1.5)&(1.5,1.866)%
          &(1,2)&(.5,1.866)&(0,1)&(1,0)&(2,0.7)}
        ,(1.932,1.0);(2.5,1.4)**\crv{(1.8,1.07)&(1.7,1.09)%
          &(1.7,1.15)&(1.93,1.35)} ,(4.55,.64)*{\bullet},(0.45,.64)*{\bullet}
        ,(4.55,.64);(5.05,0.14)**\dir{-}
        ,(0.45,.64);(-.05,0.14)**\dir{-}
        ,(4.55,1.65)*{\bullet},(0.45,1.65)*{\bullet}
        ,(3.064,1)*{\bullet},(1.938,1)*{\bullet},(1.15,1.2)*{z}
        ,(4.55,1.65);(5.05,2.15)**\dir{-}
        ,(0.45,1.65);(-.05,2.15)**\dir{-} ,(4.9,0)*{\scriptstyle{k}}
        ,(0.1,0)*{\scriptstyle{i}} ,(4.9,2.29)*{\scriptstyle{k}}
        ,(0.1,2.29)*{\scriptstyle{i}} ,(5.1,0.4)*{\scriptstyle{l}}
        ,(-0.1,0.4)*{\scriptstyle{j}} ,(5.1,1.89)*{\scriptstyle{l}}
        ,(-0.1,1.89)*{\scriptstyle{j}} }\endxyc} = -{s_1}^6
\cdot \frac{1}
    {z (z+\Lambda_i)^2 (z+\Lambda_j) (z+\Lambda_k)^2 (z+\Lambda_l)}
    \allowbreak 
    \\
    &\qquad= -{s_1}^6 \cdot 1/z^7 +{s_1}^6 (2{\Lambda_i} +
{\Lambda_j} + 2{\Lambda_k} + {\Lambda_l}) \cdot 1/z^8 +
\cdots
    \allowbreak \\
   &\qquad= 
    -{s_1}^6 \left(
      {\xyc\rgvertex[\!]4\fence2{i}{j}%
        \fence3{j}{i}\endxyc} 
      \negquad
      {\xyc\rgvertex[\!]4\fence1{k}{l}%
        \fence4{l}{k}\endxyc} 
    \right) \frac{1}{z^7}+\allowbreak \\
&\qquad\qquad
    +{s_1}^6 \left[\left(
      {\xyc\rgvertex[\varphi]4\fence2{i}{j}%
        \fence3{j}{i}\cilia{1}{2}\endxyc} 
      \negqquad
      {\xyc\rgvertex[\!]4\fence1{k}{l}%
        \fence4{l}{k}\endxyc} 
    \right) 
    + \left(
      {\xyc\rgvertex[\!]4\fence2{i}{j}%
        \fence3{j}{i}\endxyc} 
      \negqquad
      {\xyc\rgvertex[\varphi]4\fence1{k}{l}%
        \fence4{l}{k}\cilia{1}{4}\endxyc} 
    \right) \right] \cdot \frac{1}{z^8}+\cdots
   ,
\end{align*}
where
\begin{math}
    \varphi(\theta_1,\theta_2) = 2\theta_1 + \theta_2.
\end{math}
Note that, in contrast with \prettyref{xmp:hole-d}, the polynomial
$\varphi$ is not cyclically invariant.
\end{example}

\begin{example}\label{xmp:hole-f}
  \begin{align*}
    &{\xyc/r18pt/:0,
      *!C\xybox{
        ,(3.1,1.4);(3.05,1.2)**\crv{(3.13,1.5)&(3.5,1.866)%
          &(4,2)&(4.5,1.866)&(5,1)&(4,0)&(3,0.7)}
        ,(3.068,1.0);(2.5,1.4)**\crv{(3.2,1.07)&(3.3,1.09)%
          &(3.3,1.15)&(3.07,1.35)}
        ,(1.9,1.4);(1.95,1.2)**\crv{(1.87,1.5)&(1.5,1.866)%
          &(1,2)&(.5,1.866)&(0,1)&(1,0)&(2,0.7)}
        ,(1.932,1.0);(2.5,1.4)**\crv{(1.8,1.07)&(1.7,1.09)%
          &(1.7,1.15)&(1.93,1.35)}%
        ,(0.25,1);(-0.4,1)**\crv{(-0.2,1)},%
        ,(0.25,1)*{\bullet},%
        ,(3.064,1.0)*{\bullet},(1.938,1.0)*{\bullet},(1.15,1.2)*{z}
        ,(-0.4,1.2)*{\scriptstyle{i}}
        ,(-0.4,0.8)*{\scriptstyle{i}} }\endxyc} = \sum_{j=1}^N\left(
      {\xyc/r18pt/:0, *!C\xybox{
          ,(3.1,1.4);(3.05,1.2)**\crv{(3.13,1.5)&(3.5,1.866)%
            &(4,2)&(4.5,1.866)&(5,1)&(4,0)&(3,0.7)}
          ,(3.068,1.0);(2.5,1.4)**\crv{(3.2,1.07)&(3.3,1.09)%
            &(3.3,1.15)&(3.07,1.35)} 
          ,(1.9,1.4);(1.95,1.2)**\crv{(1.87,1.5)&(1.5,1.866)%
            &(1,2)&(.5,1.866)&(0,1)&(1,0)&(2,0.7)}
          ,(1.932,1.0);(2.5,1.4)**\crv{(1.8,1.07)&(1.7,1.09)%
            &(1.7,1.15)&(1.93,1.35)}%
          ,(0.25,1);(-0.4,1)**\crv{(-0.2,1)},%
          ,(0.25,1)*{\bullet},%
          ,(3.064,1.0)*{\bullet},(1.938,1.0)*{\bullet},(1.15,1.2)*{z}%
          ,(-0.4,1.2)*{\scriptstyle{i}}
          ,(-0.4,0.8)*{\scriptstyle{i}},(4.7,0.4)*{\scriptstyle{j}}
        }\endxyc}\right)
    \allowbreak \\
    &\qquad = -\U {s_1}^3\sum_{j=1}^N \cdot \frac{1}{z (z+\Lambda_i)^2 (z+\Lambda_j)}
    \allowbreak \\
    &\qquad = \U {s_1}^3 \sum_{j=1}^N \bigl( -1/z^4 +( 2\Lambda_i+\Lambda_j) \cdot 1/z^5
    + \cdots \bigr)
    \allowbreak \\
    &\qquad =-\U {s_1}^3\tr\Lambda^0\cdot 1/z^4+\U {s_1}^3\bigl(\tr\Lambda+
    2(\tr\Lambda^0)\Lambda_i\bigr) \cdot 1/z^5+\cdots
    \allowbreak \\
    &\qquad = -\U{s_1}^3\tr\Lambda^0 \Bigl(
    {\xyc\rgvertex[\!]6\fence4{i}{i}\endxyc} \negqquad\Bigr)
    \frac{1}{z^4} +
    \allowbreak \\
    & \qquad\qquad +\U{s_1}^3\Biggl[ \tr\Lambda \Bigl(
    {\xyc\rgvertex[\!]6\fence4{i}{i}\endxyc} \negqquad\Bigr) + 2\tr\Lambda^0
    \Bigl( {\xyc\rgvertex[\varphi]6\fence4{i}{i}\endxyc} \negqquad\Bigr)
    \Biggr] \frac{1}{z^5}+\cdots,
  \end{align*}
where
\begin{math}
  \varphi(\theta_1) = \theta_1.
\end{math}
Note that in this last example traces of positive powers of \(\Lambda\)
appear at the right-hand side.
\end{example}

All ribbon graphs at the right-hand side in the previous examples
are of a peculiar kind, namely, they are disjoint union of special
vertices. 
\begin{definition}
  A cluster of special vertices $\Xi$ is a ribbon graph of the form
  \begin{equation*}
    \Xi = {\sf v}_{n_1}^{\varphi_{1}} {\textstyle\coprod} \cdots
    {\textstyle\coprod} {\sf v}_{n_{k}}^{\varphi_{k}} ,
  \end{equation*}
  where the symbol $\coprod$ denotes disjoint union.  
  The valence of a cluster \(\Xi\) is the sum of the valences of its
  vertices; it is denoted by \(\text{\rm val}(\Xi)\).
  The degree of a cluster $\Xi$ is the sum of degrees of the polynomials
  decorating its vertices; denote it by $\deg\Xi$.
\end{definition}

With these notations, examples \ref{xmp:hole-a}--\ref{xmp:hole-f} show
that, for any hole type $\Gamma$, we have
\begin{equation*}\label{eq:to-be-solved-3}
  \Coeff_z^{-k} \frac{Z_{z\oplus\Lambda,s_*}^{}(\Gamma)}
    {\card{\Aut\Gamma}} 
    = \sum_{\Xi \in \X^k_\Gamma} Q^k_\Xi(s_*, \tr \Lambda^*) 
    \frac{Z_{\Lambda,s_*}^{}(\Xi)}{\card{\Aut\Xi}}\,,
\end{equation*}
where $\X^k_\Gamma$ is a suitable set of clusters of special vertices,
and $Q^k_\Xi \in \setC[s_0, s_1, \break s_2,\ldots; \tr \Lambda^0,
\tr \Lambda, \tr \Lambda^2, \ldots]$.

Moreover, we can assume that polynomials decorating special vertices
of $\Xi$ are cyclic; indeed, for any polynomial $\varphi(\theta_1,
\ldots, \theta_n)$, the cyclic polynomial
\begin{equation}
  \label{eq:to-cyclic}
  \bar \varphi(\theta_1, \ldots, \theta_n) :=  
  \sum_{\sigma \in \setZ/n\setZ} \varphi(\theta_{\sigma(1)}, 
  \ldots, \theta_{\sigma(n)})
\end{equation}
satisfies
\begin{equation}
  \label{eq:to-cyclic-g}
  \gint{{\xyc\rgvertex[\varphi]{7}\loose1\loose2\loose3\loose4\missing5\loose7%
      \cilia{1}{2}\endxyc}}
  =\gint{{\xyc\rgvertex[\bar{\varphi}]{7}\loose1%
      \loose2\loose3\loose4\missing5\loose7%
      \endxyc}}
\end{equation}
The above argument can be straightforwardly adapted to clusters
made up by several vertices.

Finally, up to splitting polynomials $\bar \varphi$ into homogeneous
components, we can further assume that polynomials decorating each
$\Xi$ are homogeneous.

The arguments used in this section lead to the following proposition,
which summarizes the way Laurent coefficients $\Coeff_z^{-k}$ transform
$z$-hole types into clusters of special vertices.
\begin{proposition}
  \label{prop:holes-to-vclusters}
  For any hole type $\Gamma$ with only ordinary vertices and each
  \(k\in\setN\) there exist:
  \begin{enumerate}
  \item a set $\X^k_\Gamma$ of clusters of vertices decorated by
    homogeneous cyclic polynomials $\varphi \in \setC[\theta_1,
    \theta_2, \ldots]$,
  \item polynomials $Q^k_\Xi \in \setC[s_0, s_1, s_2, \ldots;
    \tr\Lambda^0,\tr\Lambda,\tr\Lambda^2,\dots]$,
  \end{enumerate}
  such that:
  \begin{equation}\label{eq:coeff-to-vcluster}
    \Coeff_z^{-k}
\langle\!\langle\Gamma\rangle\!\rangle^{[1]}_{z\oplus\Lambda,s_*} =
\sum_{\Xi \in \X^k_\Gamma} Q^k_\Xi(s_*, \tr \Lambda^*) \cdot
\langle\!\langle\Xi\rangle\!\rangle_{\Lambda,s_*}.
  \end{equation}
\end{proposition}

A more accurate description of the polynomials
$Q^k_\Xi(s_*,\tr\Lambda^*)$
will be useful in the sequel of this paper.
\begin{proposition}
  \label{prop:bounds} 
  If \(\Gamma\) is a \(z\)-hole type with only ordinary vertices, then,
  with the notations of \prettyref{prop:holes-to-vclusters} above, the
  polynomials \(Q^k_\Xi(s_*,\tr\Lambda^*)\) have the form
  \begin{equation}\label{eq:valence}
    Q^k_{\Xi}(s_*,\tr\Lambda^*)= {s_*}^{m_*} \cdot q^{k}_{\Xi}(\tr\Lambda^*)\,,
  \end{equation}
  where $m_i$ is the number of $(2i+1)$-valent vertices in $\Gamma$.
  Moreover, the following inequalities hold:
  \begin{enumerate}
  \item $\text{\rm val}(\Xi) \leq \sum_{i=1}^\infty (2i-1)m_i$;\label{item:valXi}
  \item \(\deg\Xi +\deg_{s_*} Q^k_{\Xi} \leq 2k\);\label{item:degp}
  \item \(\deg\Xi +\deg^{}_{s_1}Q^k_{\Xi}\leq k\);\label{item:degptwo}
  \item if equality holds in \textsl{\ref{item:degptwo})}, then \(\Xi\)
    consists of a single vertex of valence
    \((k-\deg\Xi)\).\label{item:degpthree}
  \end{enumerate}
\end{proposition}
\begin{proof} 
  Let $m_i$ be the number of $(2i+1)$-valent vertices of $\Gamma$; equation
  \eqref{eq:valence} follows by a straightforward application of the
  Feynman rules; so we only need to show bounds
  \textsl{\ref{item:valXi})}--\textsl{\ref{item:degpthree})}.
  
  The valence (i.e., the number of legs) of any cluster of vertices
  \(\Xi\) at right-hand side in \prettyref{eq:coeff-to-vcluster} is
  exactly the number of half-edges which stem from the vertices of
  \(\Gamma\) and which do not border the \(z\)-decorated hole.  If a
  half-edge of \(\Gamma\) stems from a 1-valent vertex, then it must border
  the \(z\)-decorated hole on both sides; when $i\geq1$, at most
  \((2i-1)\) half-edges stemming from a \((2i+1)\)-valent vertex may
  not border the \(z\)-decorated hole. This proves
  \textsl{\ref{item:valXi})}.

  Let \(\nu\) be the number of internal edges of \(\Gamma\). Since an
  internal edge of $\Gamma$ carries either a factor $1/z$ or a factor
  \(2/(z+\Lambda_i)=(2/z)(1-\Lambda_i/z+{\Lambda_i}^2/z^2-\cdots)\),
  then the Laurent coefficient of \(z^{-k}\) is a polynomial of degree
  at most \(k-\nu\) in the $\Lambda_i$'s.  The graph \(\Gamma\) can
  have at most \(2\nu\) vertices, and the Laurent coefficient of
  \(z^{-k}\) can be non-zero only if \(k \geq \nu\). Then
  \(\deg\Xi+\deg_{s_*} Q^k_{\Xi}\leq (k-\nu)+2\nu\leq 2k\), which is
  \textsl{\ref{item:degp})}.
  
  If \(\Gamma\) has \(m_1\) trivalent vertices, then it has at least
  \(m_1\) internal edges, so \(\deg\Xi+\deg_{s_1} Q^k_{\Xi}\leq
  (k-m_1)+m_1=k\). Then \textsl{\ref{item:degptwo})} is proven.
  
  Assume now that \(\deg\Xi + \deg_{s_1} Q^k_{\Xi} = k\). The number
  \(m_1\) of 3-valent vertices of \(\Gamma\) is \(\deg_{s_1} Q^k_{\Xi}\).
  Since \(m_1\) is at most equal to the number \(\nu\) of internal
  edges of \(\Gamma\), we have:
  \begin{equation*}
    k= \deg\Xi + \deg_{s_1} Q^k_{\Xi} \leq
    (k-\nu) + m_1  \leq  k,
  \end{equation*}
  which forces $m_1 =\nu$ and \(\deg\Xi = k-\nu\).  Since \(\deg\Xi =
  k-\nu\), no edge of \(\Gamma\) borders the \(z\)-decorated hole on both
  sides; this, together with $m_1 =\nu$, implies that all the vertices of
  \(\Gamma\) are trivalent and \(\Gamma\) must be the \(z\)-hole type
  \begin{equation*}
    {\xy
      \rghole[z]{7}\loose1\loose2\loose3\loose4\loose5\loose6\loose7
      \endxy}
    \qquad
    \text{($\nu$ legs)}
  \end{equation*}
  which is changed into a \(\nu\)-valent special vertex by the operation
  of taking a coefficient of the Laurent expansion of its amplitude
  with respect to \(1/z\) at \(z=\infty\).
\end{proof}

\subsection{Expectation values of polynomial vertices}
\label{sec:expval-poly}

The aim of this section is to find a canonical form to express
expectation values of polynomial-decorated vertices.

\begin{definition}
  Let \(\varphi \in \setC[\theta_1,\theta_2,\dots,\theta_n]\). We say
  that $\varphi$ is cyclically decomposable iff there exists $\psi \in
  \setC[\theta_1,\theta_2,\dots,\theta_n]$ such that
  \begin{equation*}
    \varphi(\theta_1,\theta_2,\dots,\theta_n) = 
    \sum_{\sigma\in\setZ/n\setZ} (\theta_{\sigma(n)} + \theta_{\sigma(1)}) 
    \psi(\theta_{\sigma(1)}, \theta_{\sigma(2)}, \dots, \theta_{\sigma(n)}).
  \end{equation*}
  
  We say that $\varphi$ is residual iff it is a degree zero polynomial or
  has the form:
  \begin{equation*}
    \varphi(\theta_1, \ldots, \theta_{2n}) = \text{const} \cdot \sum
\theta_i^{2d}.
  \end{equation*}
\end{definition}

\begin{lemma}
  \label{lemma:cyclic-decomposition}
  Every homogeneous cyclic polynomial 
$\varphi \in \setC[\theta_1, \ldots, \theta_n]$ can
  be split into a sum of a cyclically decomposable $\varphi\dec$ and a
  residual $\varphi\res$:
\begin{equation*}
  \varphi (\theta_1, \ldots, \theta_n) =   \varphi\dec (\theta_1, \ldots,
  \theta_n) +   \varphi\res (\theta_1, \ldots, \theta_n) ,
\end{equation*}
\end{lemma}
\begin{proof}
Let $d$ be the degree of $\varphi$. The statement is trivial if $d=0$, so
assume $d\geq1$. Let $I^d_n$ be the ideal in $\setC[\theta_1, \ldots,
\theta_n]$ generated by
$\theta_1+\theta_2,\theta_2+\theta_3,\dots,\theta_n+\theta_1$. If $n$ is
odd, then  $I_n^d=\setC[\theta_1, \ldots, \theta_n]$, and from the
cyclical invariance of $\varphi$ it easily follows that it is cyclically
decomposable. If $n$ is even, then $I_n^d$ is the ideal of polynomials
that vanish at the point $(1,-1,1,\dots,-1)$. Therefore, by adding to
$\varphi$ a
suitable multiple of $\sum
\theta_i^{2d}$ we get an element of $I_n^d$ which is cyclically invariant,
and so cyclically decomposable.
\end{proof}

In particular, a polynomial of positive degree in an odd number of
indeterminates is always cyclically decomposable.

\begin{definition}\label{dfn:decomposable-vs-residual}
  We say that a cluster \(\Xi\) is \emph{decomposable} if at
  least one vertex of $\Xi$ is decorated by a cyclically decomposable
  polynomial, otherwise we say that $\Xi$ is \emph{residual}.
\end{definition}
By linearity, each cluster $\Xi$ can be split into a sum $\Xi =
\Xi\dec + \Xi\res$ where \(\Xi\dec\) is decomposable and \(\Xi\res\)
is residual.

Motivation for distinguishing between decomposable and residual
clusters is given by \prettyref{prop:contraction}; to prove it, we
need first a technical result.
\begin{lemma}  
  \label{lemma:lower-degree}
  If $\Upsilon$ is a decomposable cluster, then its expectation value
  can be written as a linear combination (over $\setC[s_*,
  \tr\Lambda^*]$) of expectation values of clusters of lower
  degree:
  \begin{equation*}
    \expval{\Upsilon}^{}_{\Lambda,s_*} = \sum_{\Xi \in \X_\Upsilon}
p_\Xi^{}(s_*, \tr\Lambda^*) \expval{\Xi}^{}_{\Lambda,s_*},
  \end{equation*}
  with
  \begin{equation*}
    \deg\Xi +\deg^{}_{s_*}p^{}_{\Xi}\leq \deg\Upsilon,\quad
    \forall\Xi \in \X_\Upsilon.
  \end{equation*}
\end{lemma}
\begin{proof} 
  We first give a proof for a cluster made up of
  a single vertex.
  
  For any $\psi \in \setC[\theta_1, \ldots, \theta_n]$, let
  \(u_\psi \in \setC[\theta_1,\dots,\theta_n]\) be the polynomial
  \begin{equation*}
    u_\psi(\theta_1,\theta_2,\dots,\theta_n):=
    (\theta_n+\theta_1)\cdot\psi(\theta_1,\theta_2,\dots,\theta_n) .
  \end{equation*}
 Then $\varphi$ is cyclically decomposable iff, for some $\psi$:
  \begin{equation*}
    \varphi(\theta_1,\dots,\theta_n) = \sum_{\sigma\in\setZ/n\setZ} u_\psi(
    \theta_{\sigma(1)}, \dots, \theta_{\sigma(n)});
  \end{equation*}
  this implies the graphical identity
  \begin{equation*}
    \gint{{\xyc\rgvertex[\varphi]{7}\loose1\loose2\loose3\loose4\missing5\loose7%
        \endxyc}}
    =\gint{{\xyc\rgvertex[u_\psi]{7}%
        \loose1\loose2\loose3\loose4\missing5\loose7%
        \cilia{1}{2}\endxyc}} .
  \end{equation*} 
  By definition of expectation value of a diagram, both sides are sums
  over ribbon graphs (with distinguished sub-diagrams); for any \(\Gamma\)
  in the sum at right-hand side, the edge stemming from the vertex
  just before the ciliation (in the cyclic order of the vertex) must
  \emph{either} end at another ---distinct--- vertex \emph{or} make a
  loop. Therefore, using the definition of expectation value again,
  \begin{multline}\label{eq:expansion}
    \gint{{\xyc\rgvertex[u_\psi]{13}\loose1%
        \loose3%
        \loose5%
        \loose{10}%
        \loose{12}%
        \missing7%
        \cilia{1}{3}%
        \endxyc}}=\sum_{j=0}^{\infty}\gint{{
        \xyc
        \rgvertex[u_\psi]{13}%
        \loose1%
        \loose3%
        \loose5%
        \loose{10}%
        \loose{12}%
        \missing7%
        \cilia{1}{3}%
        ,(1.3,0);(1.6,.5)**\dir{-},(1.3,0);(1.6,-.5)**\dir{-}%
        ,(1.3,0)*{\bullet}%
        ,(1.7,0)*{\scriptstyle{\cdots}}
        ,(2.6,0)*{\bigg\}\text{$2j$ edges}}
        \endxyc}}
    \\
    +\gint{{\xyc
        \rgvertex[u_\psi]{13}%
        \join{1}{3}
        \loose5%
        \loose{10}%
        \loose{12}%
        \missing7%
        \cilia{1}{3}%
        \endxyc}}
    +\sum_{j=1}^{n-3}\gint{{
        \xyc
        \rgvertex[u_\psi]{13}%
        \loose2%
        \loose4%
        \loose6%
        \loose{10}%
        \missing{3}%
        \loose{12}%
        \missing7%
        \cilia{1}{2}%
        ,(1.5,1.5)*{j\ {\rm edges,}}
        ,(1.5,1.2)*{0<j<n-2}
        ,(.405,0);(.85,.35)**\crv{(1.2,-0.2)}%
        ,(.75,.48);(.15,.92)**\crv{(.6,.7)}
        ,(-.15,.39);(.05,.95)**\crv{(-.3,1.1)}
        \endxyc
      }}
    \\
    +\gint{{\xyc
        \rgvertex[u_\psi]{13}%
        \join{1}{12}
        \loose5%
        \loose{10}%
        \loose{3}%
        \missing7%
        \cilia{1}{3}%
        \endxyc}}
  \end{multline}
  Now, each of the terms at right-hand side of \eqref{eq:expansion}
  above, can be rewritten as the expectation
  value of a linear combination (over $\setC[s_*; \tr \Lambda^*]$) of
  clusters; indeed, one can directly compute:
  \begin{align}
    \label{eq:C1}\tag{C1}
    \xyc
    \rgvertex[u_\psi]{13}%
    \loose1%
    \loose3%
    \loose5%
    \loose{10}%
    \loose{12}%
    \missing7%
    \cilia{1}{3}%
    ,(1.3,0);(1.6,.5)**\dir{-},(1.3,0);(1.6,-.5)**\dir{-}%
    ,(1.3,0)*{\bullet}
    ,(1.7,0)*{\scriptstyle{\cdots}}
    ,(2.6,0)*{\bigg\}\text{$2j$ edges}}
    \endxyc
    &= 
    \U\Biggl(-\frac{1}{2}\Biggr)^{j-1}s_j \cdot
    \underbrace{
      \xyc
      \rgvertex[\psi]{13}%
      \missing1%
      \loose2%
      \loose{13}%
      \loose3%
      \loose5%
      \loose{10}%
      \loose{12}%
      \missing7%
      \cilia{2}{3}
      \endxyc
    }_{\text{$(n+2j-1)$-valent}},
    \\ \label{eq:C2}\tag{C2}
    \xyc
    \rgvertex[u_\psi]{13}%
    \join{1}{3}
    \loose5%
    \loose{10}%
    \loose{12}%
    \missing7%
    \cilia{1}{3}%
    \endxyc
    &= 
    2\sum_{h=0}^k\tr\Lambda^h\cdot
    \underbrace{
      \xyc
      \rgvertex[\psi_h]{13}%
      \loose5%
      \loose{10}%
      \loose{12}%
      \missing7%
      ,(.5,.5)*{\scriptstyle{\bigstar}}
      \endxyc
    }_{\text{$(n-2)$-valent}}
    ,
    \\
    \label{eq:C3}\tag{C3}
    {\xyc
      \rgvertex[u_\psi]{13}%
      \loose2%
      \loose4%
      \loose6%
      \loose{10}%
      \missing{3}%
      \loose{12}%
      \missing7%
      \cilia{1}{2}%
      ,(1.5,1.5)*{j\ {\rm edges,}}
      ,(1.5,1.2)*{0<j<n-2}
      ,(.405,0);(.85,.35)**\crv{(1.2,-0.2)}%
      ,(.75,.48);(.15,.92)**\crv{(.6,.7)}
      ,(-.15,.39);(.05,.95)**\crv{(-.3,1.1)}
      \endxyc}
    &= 
    2\sum_{h=0}^k
    \underbrace{
      \xyc
      \rgvertex[\phi_{h}']{13}%
      \loose2%
      \loose4%
      \missing3%
      \ciliamedio{1}{2}%
      \endxyc
    }_{\text{$j$-valent}}
    \ \underbrace{
      \xyc
      \rgvertex[\phi_{h}'']{13}%
      \loose6%
      \loose{12}\loose{10}%
      \missing7%
      \ciliamedio{1}{2} 
      \endxyc
    }_{\text{$(n-j-2)$-valent}}
    ,
  \end{align}
  \begin{align}
    \label{eq:C4}\tag{C4}
    \xyc
    \rgvertex[u_\psi]{13}%
    \join{1}{12}
    \loose5%
    \loose{10}%
    \loose{3}%
    \missing7%
    \cilia{1}{3}%
    \endxyc
    &= 
    2\sum_{h=0}^k\tr\Lambda^h\cdot
    \underbrace{
      \xyc
      \rgvertex[\eta_h]{13}%
      \loose5%
      \loose{10}%
      \loose{3}%
      \missing7%
      ,(.5,-.5)*{\scriptstyle{\bigstar}}
      \endxyc
    }_{\text{$(n-2)$-valent}}
    , 
  \end{align}
  for polynomials $\psi_h$, $\phi_{h}'$, $\phi_{h}''$,
  $\eta_h$
  defined by:
  \begin{gather*}
    \sum_{h=0}^k \theta_1^h \psi_h(\theta_2, \dots, \theta_{n-1}) =
    \psi(\theta_1, \dots, \theta_{n-1}, \theta_2) ,
    \\
    \sum_{h=0}^k \phi_{h}'(\theta_1, \dots, \theta_j) \cdot
    \phi_{h}''(\theta_{j+2}, \dots, \theta_{n-1}) 
    \phantom{XXXXXXXXXX}
    \\
    \phantom{XXXXXXXXXX}
    = \psi(\theta_1,
    \dots, \theta_j, \theta_1, \theta_{j+2}, \dots, \theta_{n-1},
    \theta_{j+2}) ,
    \\
    \sum_{h=0}^k \theta_n^h \eta_h(\theta_1, \dots, \theta_{n-2}) =
    \psi(\theta_1, \dots, \theta_{n-2}, \theta_1, \theta_n).
  \end{gather*}
  
  The general case of clusters made up of more than $1$ vertex
  is done by picking a vertex out of the cluster and applying the
  above procedure to it. A new combination of vertices may appear, which
  is not listed in equations \eqref{eq:C1}--\eqref{eq:C4} above;
  namely, that the ciliated edge connects the chosen vertex to another
  one in the same cluster. Direct computation again gives:
  \begin{equation*}
    \label{eq:C1'}\tag{C1'}
    \xyc
    \rgvertex[u_\psi]{13}%
    \loose1%
    \loose3%
    \loose5%
    \loose{10}%
    \loose{12}%
    \missing7%
    \cilia{1}{3}%
    ,(1.6,0)*{\ \zeta\ }*\cir{};(2.15,.8)**\dir{-}%
    ,(1.6,0)*{\ \zeta\ }*\cir{};(2.15,-.8)**\dir{-}%
    ,(2.1,0)*{\scriptstyle{\cdots}}
    ,(3.2,0)*{\Bigg\}\text{$j$ edges}}
    \endxyc
    = 2\cdot\underbrace{
      \xyc
      \rgvertex[\psi*\zeta]{13}%
      \loose2%
      \loose{13}%
      \loose3%
      \loose5%
      \loose{10}%
      \loose{12}%
      \ciliamedio{2}{3},(-1,.2)*{\cdots},(1,0)*{\cdots}
      \endxyc
    }_{\text{$(n+j-1)$-valent}}
    ,
  \end{equation*}
  for a polynomial $\psi*\zeta$ given by
  \begin{equation*}
    (\psi*\zeta)(\theta_1,\dots,\theta_{n+j})= 
    \psi(\theta_1,\dots,\theta_n)\cdot\zeta(\theta_n,\theta_{n+1},
    \dots,\theta_{n+j},\theta_1).
  \end{equation*}
  This proves the claim.
\end{proof}

By repeatedly applying the edge-contraction procedure from 
\prettyref{lemma:lower-degree} to the right hand side of equation
\prettyref{eq:coeff-to-vcluster}, and by inequalities described in
\prettyref{prop:bounds}, one can prove the following.
\begin{proposition}
  \label{prop:contraction}
For any hole type \(\Gamma\) with only ordinary vertices,
and any positive integer \(k\), 
\begin{equation*}\label{eq:formulacoeff}
  \Coeff_z^{-k}\langle\!\langle\Gamma\rangle\!\rangle^{[1]}_{z\oplus\Lambda}
  = \sum_{\norm{m_*}\leq 2k} \sum_{\Xi\in{\mathcal{X}}_{m_*,\Gamma}}  
  s_*^{m_*} q_{m_*,\Xi}^{}(\tr\Lambda^*) \langle\!\langle\Xi\rangle\!\rangle^{}_{\Lambda}
\end{equation*}
for suitable residual clusters \(\Xi\) and polynomials
\(q_{m*,\Xi} \in \setC[\tr\Lambda^0, \tr \Lambda, \break \tr
\Lambda^2, \ldots]\).  Moreover, for any cluster in \({\mathcal
  X}_{m_*,\Gamma}\), the following inequalities hold:
\begin{enumerate}
\item\label{item:i1} \(\text{\rm val}({\Xi})\leq \sum_i (2i-1)m_i\); 
\item \(\deg\Xi +\abs{m_*} \leq 2k\);
\item\label{item:i3} \(\deg\Xi +m_1\leq k\) and if \(\deg\Xi +m_1= k\)
  then \(\Xi\) consists of a single special vertex of valence
  \((k-\deg\Xi\)).
\end{enumerate}  
\end{proposition}

\subsection{Proof of the Main Theorem}

To conclude the proof of the main result of this paper, we need to
introduce an algebra of formal differential operators in the variables
$s_*$.  For any polyindex \(m_*=(m_0,m_1,\dots,m_l,0,0,\dots)\) set:
\begin{equation*}
\abs{m_*} := \sum_{i=0}^\infty m_i,\quad
\norm{m_*}_-:=\sum_{i=1}^\infty(2i-1)m_i,\quad
\norm{m_*}_+:=\sum_{i=0}^\infty(2i+1)m_i.
\end{equation*}
\begin{definition}
  \label{dfn:formal-diff-operator}
  A formal triangular differential operator in the variables
$s_*$ is a
formal
  series
  \begin{equation*}
    D(s_*,\del / \del{s_*}) = \sum_{\norm{n_*}_+ \leq \norm{m_*}_-}
    a_{m_*,n_*} s_*^{m_*} \frac{\del^{\abs{n_*}}}{\del s_*^{n_*}},
    \qquad a_{m_*, n_*} \in \setC,
  \end{equation*}
  of bounded degree in $s_*$.
 \end{definition}
Formal triangular differential operators in the variables $s_*$ form a
(non commutative) algebra $\setC \langle\!\langle s_*,
\del / \del{s_*} \rangle\!\rangle$.

\begin{theorem}\label{theorem:main}
  For any $k\geq0$ there exists a formal triangular differential operator
\begin{equation*}
    D_{k} = c_k s_1^{2k+1} {\del / \del s_{k}} 
    + \text{lower $s_1$-degree terms},\quad c_k\in\setC.
  \end{equation*}
 such that
  \begin{equation*}
    \frac{\del Z(s_*;t_*)}{\del t_k}
    =D_{k}(s_*,\del / \del{s_*})Z(s_*;t_*).
  \end{equation*}
\end{theorem}
\begin{proof}
  By equation \prettyref{eq:dZst-twotwo},
\begin{equation}
  \label{eq:derivatives-partition}
  \left.\frac{\del Z(s_*;t_*)}{\del t_k}\right|_{t_*(\Lambda)}
  = -\frac{1}{(2k-1)!!}\sum_{\Gamma\in{\mathcal
S}}\Coeff_z^{-(2k+1)}\langle
\!\langle\Gamma\rangle\!\rangle^{[1]}_{z\oplus\Lambda,s_*}\,,
\end{equation} 
where \(\mathcal{S}\) denotes the set of all the $z$-hole types.  

By the Feynman rules for the $(N+1)$-dimensional 't~Hooft-Kontsevich
model, each edge with one or both sides decorated by the variable
$z$ corresponds to a factor of order $O(z^{-1})$ as $z\to\infty$.
This implies that $Z_{z\oplus\Lambda;s_*}(\Gamma)=O(z^{-k})$ if the
hole type $\Gamma$ is such that the $z$-decorated hole is bounded by
$k$ edges, so that
\(\Coeff_z^{-(2k+1)}\bigl(Z_{z\oplus\Lambda;s_*}(\Gamma)\bigr)=0\)
if more than \(2k+1\) edges border the \(z\)-decorated hole. Equation
\prettyref{eq:derivatives-partition} is therefore equivalent to
\begin{equation}\label{eq:derivatives-partition2}
\left.\frac{\del Z(s_*;t_*)}{\del t_k}\right|_{t_*(\Lambda)}
=\sum_{h\leq 2k+1}\sum_{\Gamma\in{\mathcal
S}_h}\Coeff_z^{-(2k+1)}\langle\!\langle
\Gamma\rangle\!\rangle^{[1]}_{z\oplus\Lambda;s_*}\,,
\end{equation} 
where \({\mathcal S}_h\) denotes the set of hole-types whose
\(z\)-decorated hole is bounded by exactly \(h\) edges.

By applying \prettyref{prop:contraction} to the right-hand side of
equation \prettyref{eq:derivatives-partition2}, we find:
\begin{equation}\label{eq:almost-last-equation}
 \left.\frac{\del Z(s_*;t_*)}{\del t_k}\right|_{t_*(\Lambda)}
      = \sum_{\norm{m_*}\leq 4k+2}\Biggl(\sum_{\Xi \in
\X[m_*]} s_*^{m_*}q_{m_*,\Xi}^{}(\tr
\Lambda^*)\langle\!\langle{\Xi}\rangle\!\rangle_{\Lambda,s_*}^{}
\Biggr),
\end{equation}
for suitable residual clusters \(\Xi\) and polynomials
$q^{}_{m_*\Xi}\in \setC[\tr \Lambda^0, \tr \Lambda^1, \break \tr
\Lambda^2, \ldots]$. The behavior of the left hand side imposes
strict constraints both on $q_{m_*, \Xi}$ and the clusters
$\Xi$.

\textsl{1)} The polynomials $q_{m_*, \Xi}(\tr \Lambda^*)$ are constant
with
respect to $\Lambda$; indeed, since $Z(s_*;t_*)$ is a formal power series
in the variables $t_*$, the left-hand side of
\prettyref{eq:almost-last-equation} is a function of the traces of
negative powers of \(\Lambda\) only; thus, terms in the right-hand side
containing traces of positive powers of
\(\Lambda\) must cancel out.  Therefore, $q_{m_*, \Xi} (\tr \Lambda^*) =
r_{m_*, \Xi}
\in \setC$, so that
\begin{equation}\label{eq:last-equation}
  \left.
    \frac{\del Z(s_*;t_*)}{\del t_k}\right|_{t_*(\Lambda)}
  = \sum_{\norm{m_*}\leq 4k+2}\Biggl(\sum_{\Xi \in \X[m_*]}
r_{m_*,\Xi}^{}\,s_*^{m_*}\langle\!\langle{\Xi}\rangle\!\rangle_{\Lambda,s_*}^{}
  \Biggr).
\end{equation}

\textsl{2)} All the vertices appearing in the clusters on the right
hand side of equation \prettyref{eq:last-equation} are odd-valent.
Indeed, since formula \prettyref{eq:last-equation} holds for every
\(N\), a fortiori it holds for \(N=2\); both sides of
\prettyref{eq:last-equation} are real analytic for positive real
\(\Lambda_1\),\(\Lambda_2\); their analytic prolongations coincide on
the connected region
$U=\setC^2\setminus\{\Lambda_1=0;\Lambda_2=0;\Lambda_1+\Lambda_2=0\}$.
For real positive $\varepsilon$, set
  \begin{equation*}
    \Lambda^{}_\varepsilon(\lambda) :=
    \begin{pmatrix}
      -\lambda + 2\varepsilon  &  0
      \\
      0                          &  \lambda
    \end{pmatrix} .
  \end{equation*}
  For any \(|\lambda|>2\varepsilon\), the diagonal matrix
\(\Lambda_\varepsilon^{}(\lambda)\) lies in \(U\)
  and we can consider \prettyref{eq:last-equation} at
\(\Lambda=\Lambda_\varepsilon^{}(\lambda)\).

Since $t_*\bigl(\Lambda_\varepsilon(\lambda)\bigr)\to 0$ as
$\lambda\to+\infty$, the left-hand side of
\prettyref{eq:last-equation} has a finite limit, independent of
\(\varepsilon\), for \(\lambda\to+\infty\); in particular, there exist
some formal power series \(\chi_k(s_*)\in\setC[[s_*]]\), such that
\begin{equation}\label{eq:to-be-contrasted}
\lim_{\varepsilon\to 0}\lim_{\lambda\to+\infty}
\sum_{\norm{m_*}\leq 4k+2}\Biggl(\sum_{\Xi \in
\X[m_*]}
r_{m_*,\Xi}^{}\,s_*^{m_*}\langle\!\langle{\Xi}\rangle\!\rangle_
{\Lambda_\varepsilon^{}(\lambda),s_*}^{}\Biggr)=\chi_k(s_*)\,.
\end{equation}
Assume an even-valent vertex \({\sf v}_{2n}^\varphi\) appears in a
cluster \(\Xi_0\) on the right-hand side of
\prettyref{eq:last-equation}, and let \(d\) be the degree of the
polynomial \(\varphi(\theta_1\dots\theta_{2n})\). Since all
clusters on the right-hand side of \prettyref{eq:last-equation} are
residual, the degree \(d\) is even and
\begin{equation*}
  \varphi(\theta_1\dots\theta_{2n})= \text{const} \cdot
\sum_{i=1}^{2n}{\theta_i}^d \,.
\end{equation*}
According to the Feynman rules for the 't~Hooft-Kontsevich model, the
expectation value
\(\langle\!\langle\Xi_0\rangle\!\rangle^{}_{\Lambda^{}_\varepsilon(\lambda)
,s_*}\)
expands into a sum over ribbon graphs whose holes are colored with the two
colors \(1\),
\(2\).  The edges of such a graph fall within one of these kinds:
\begin{enumerate}
\item\label{firstype} both sides of the edge are decorated by the
  color \(1\): this edge brings a factor
  \(-1/(\lambda-2\varepsilon)\);
\item\label{secondtype} both sides of the edge are decorated by the
  color \(2\): this edge brings a factor \(1/\lambda\);
\item\label{thirdtype} one side of the edge is decorated by the color
  \(1\) and the other by the color \(2\): this edge brings a factor
  \(1/\varepsilon\).
\end{enumerate}

Since \({\sf v}_{2n}^{\varphi}\) is an even-valent vertex, in the
expansion
of
\(\langle\!\langle\Xi_0\rangle\!\rangle^{}_{\Lambda^{}_\varepsilon(\lambda)
,s_*}\)
into ribbon graphs with holes decorated by the indices \(1\) and \(2\), we
find terms with a connected component having only edges of the third type,
e.g.,
\begin{equation*}\label{eq:colouredvertex}
    {\xy\rgvertex[\varphi]{8}%
      \join12\join34%
      \join56\join78%
      ,(.5,.5)*{1}%
      ,(0,.73)*{2},(-.5,.5)*{1},(-.73,0)*{2}%
      ,(-.5,-.5)*{1},(0,-0.73)*{2},(.5,-.5)*{1},(.73,0)*{2}\endxy}
  \end{equation*}
which evaluate to
\begin{equation}
  \label{eq:divergent-term}
  (1/\varepsilon^n) \cdot \varphi(-\lambda + \varepsilon,
  \lambda, \ldots, -\lambda+ \varepsilon,
  \lambda) .
\end{equation}
If \(d>0\), this diverges as \(\lambda\to+\infty\). If \(d=0\), then in the limit
\(\lambda\to+\infty\), the term \prettyref{eq:divergent-term} has a polar behavior
as \(\varepsilon\to 0\). In either case we would have a divergent behavior
contradicting equation \prettyref{eq:to-be-contrasted}. Therefore, no
even-valent vertices can appear in the clusters on the right-hand side
of \prettyref{eq:last-equation}.

Since a residual odd-valent vertex must have degree zero, and any
cluster made up of degree zero odd-valent vertices is of the
form 
\begin{equation*}
    {{\sf v}_1}^{\coprod n_0} {\textstyle\coprod} \cdots
    {\textstyle\coprod} {{\sf v}_{2l+1}}^{\coprod n_l}
  \end{equation*}
  for some polyindex \(n_*\), we have finally proven: in the large
  \(N\) limit,
\begin{equation*}
  \left.\frac{\del Z(s_*;t_*)}{\del t_k}\right\rvert_{t_*(\Lambda)}
  = \sum_{\substack{\norm{m_*}\leq 4k+2 \\ \norm{n_*}_+ \leq
\norm{m_*}_-}}
  r_{m_*,n_*}^{}\,s_*^{m_*}\expval{
    {{\sf v}_1}^{\coprod n_0} {\textstyle\coprod} \cdots
    {\textstyle\coprod} {{\sf v}_{2l+1}}^{\coprod n_l}
  }_{\Lambda,s_*}^{}.
\end{equation*}
By equation \prettyref{eq:dnZ}, expectation values of
clusters of degree zero odd-valent special
vertices can be expressed as derivatives of the partition function of the
't~Hooft-Kontsevich model with respect to the $s_*$ variables. Namely,
\begin{equation*}
 \expval{
    {{\sf v}_1}^{\coprod n_0} {\textstyle\coprod} \cdots
    {\textstyle\coprod} {{\sf v}_{2l+1}}^{\coprod n_l}
  }_{\Lambda,s_*}^{}
  =
  \frac{\U^{|n_*|}(-2)^{\sum_j jn_j}} {n_0! \cdots n_l!}
  \frac{\del^{|n_*|}} {{\del s_0}^{n_0} \cdots {\del s_l}^{n_l}}
  \expval{\emptyset}_{\Lambda,s_*}^{}.
\end{equation*}
Therefore, there exist a formal triangular differential operator
$D_k(s_*,\del / \del{s_*})$ such that, in the large \(N\) limit,
\begin{equation*}
  \left.\frac{\del Z(s_*;t_*)}{\del t_k}\right\rvert_{t_*(\Lambda)}
  = D_k(s_*,\del / \del{s_*})
Z(s_*;t_*)\biggr\rvert_{t_*(\Lambda)}
\end{equation*}
Moreover, bounds {\it \ref{item:i1})--\ref{item:i3})} in
\prettyref{prop:contraction} dictate that $D_k$ has the form
\begin{equation*}
  D_k = c_k \cdot {s_1}^{2k+1} \partial / \partial s_k + \text{ lower
$s_1$-degree terms.}
\end{equation*}
Since, in the large $N$ limit, the $t_*(\Lambda)$ become independent
coordinates, the statement follows.
\end{proof}
Let now $\setC\langle\del/\del t_*\rangle$ be the free non-commutative
algebra generated by the $\del/\del t_k$. It acts on the space of
formal power series in the variables $t_*$ through its abelianization
$\setC[\del/\del t_*]$.

\begin{theorem}\label{theorem:2}
The map
\begin{equation*}
  \frac{\del}{\del t_k}\mapsto D_k(s_*,\del/\del s_*)
\end{equation*} 
induces an algebra homomorphism
\begin{align*}
D\colon \setC\langle \del/\del t_*\rangle&\to
\setC\langle\!\langle s_*,\del/\del
s_*\rangle\!\rangle\\
P&\mapsto D_P
\end{align*}
such that:
  \begin{equation}\label{eq:t2}
    P(\del / \del{t_*})Z(s_*;t_*)
    =D_{P}(s_*,\del /\del{s_*})Z(s_*;t_*).
  \end{equation}
 \end{theorem}
\begin{proof}  
  The statement immediately follows by the fact that, when regarded as
  differential operators acting on the space of formal power series in
  the variables $t_*$ and $s_*$, the elements $\del/\del t_k$ commute
  among themselves and with the operators $D_{P}(s_*,\del / \del{s_*})$. For
  instance, to prove that $(\del_{t_i}\cdot
  \del_{t_j})Z(s_*;t_*)=D_{\del_{t_i}\cdot
    \del_{t_j}}(s_*,\del_{s_*})Z(s_*;t_*)$, one computes:
\begin{align*}
(\del_{t_i}\cdot \del_{t_j})Z(s_*;t_*)&=
\frac{\del^2 Z(s_*;t_*)}{\del t_i\del t_j}=\frac{\del}{\del t_j}\biggl(
\frac{\del}{\del t_i}Z(s_*;t_*)\biggr)\\
&=\frac{\del}{\del
t_j}D_i(s_*,\del_{s_*})Z(s_*;t_*)\\
&=D_i(s_*,\del_{s_*})\frac{\del}{\del t_j}Z(s_*;t_*)\\
&=D_i(s_*,\del_{s_*})D_j(s_*,\del_{s_*})Z(s_*;t_*)\\
&=D_{\del_{t_i}\cdot \del_{t_j}}(s_*,\del_{s_*})Z(s_*;t_*).
\end{align*}
The proof for a higher order monomial in the $\del/\del t_*$ goes
along the same lines.
\end{proof}

\section{Examples and Applications}
\label{sec:final}

Let $s^\circ_0, \dots, s^\circ_r$ be complex constants, and set
$s^\circ_*=(s^\circ_0,s^\circ_1,\dots,s^\circ_r,0,0,\dots)$.  There
is a well-defined evaluation map
\begin{equation*}
  \ev_{s^\circ_*} : \setC\langle\!\langle s_*; \del / \del s_* \rangle\!\rangle \to \setC [\del / \del s_*],
\end{equation*}
which is linear but \emph{not} an algebra homomorphism.

\begin{corollary}\label{cor:dfiz+der}
  For any $s^\circ_* = (s^\circ_0, \dots, s^\circ_r, 0, 0, \dots)$,
  there exists a linear map
  \begin{equation}
    \label{eq:dfn-Q}
    \begin{split}
      Q^{s^\circ_*} : \setC[\del / \del t_*]
      &\to \setC [\del / \del s_*],
      \\
      P &\mapsto Q^{s^\circ_*}_P,
    \end{split}
  \end{equation}
  such that
  \begin{equation}
    \label{eq:cor1}
    P (\del / \del t_*)Z(s_*^\circ;t_*) = \ev_{s^\circ_*} \left[ 
      Q^{s^\circ_*}_P(\del /\del s_*) Z(s_*;t_*) 
    \right]
\end{equation}
\end{corollary}
\begin{proof}
  The basis $\{ \del^{\abs{m_*}} t / {\del t_0}^{m_0} {\del t_1}^{m_1}
  \cdots {\del t_q}^{m_q} \}$ defines a linear section $\varsigma$ to
  the projection $\setC \langle \del / \del t_* \rangle \to \setC[
  \del / \del t_*]$; take $D: \setC \langle \del / \del t_* \rangle
  \to \setC \langle\!\langle s_*; \del / \del s_* \rangle\!\rangle$ as
  in \prettyref{theorem:2} and set:
  \begin{equation*}
    Q^{s^\circ_*} = \ev_{s^\circ_*} \circ D \circ \varsigma.
  \end{equation*}
  Equation \eqref{eq:cor1} now follows from \prettyref{eq:t2}.
\end{proof}

Recall that $Z(t_*)$ is the partition function for intersection
numbers on the moduli spaces of stable curves; it is a special case of
$Z(s_*; t_*)$ when $s_* = (0, 1, 0, 0, \dots)$.  Thus, in particular,
we get the following.
\begin{corollary}[DFIZ Theorem]\label{cor:dfiz-der}
  There exists a linear isomorphism 
  \begin{equation*}
    Q: \setC [\del / \del t_*] \to \setC [\del / \del s_*]
  \end{equation*}
  such that 
  \begin{equation*}
    P(\del / \del t_*) Z(t_*) = \ev_{(0,1,0,0,\dots)} \left[
      Q_P(\del / \del s_*) Z(s_*; t_*)
      \right].
  \end{equation*}
\end{corollary}
\begin{proof}
  We just need to prove that the map \(Q := Q^{(0,1,0,0,\dots)}\) is a
  linear isomorphism. Indeed, from
\begin{equation*}
  D_{\del/ \del t_k} = c_k {s_1}^{2k+1} {\del / \del s_{k}} +
  \text{lower $s_1$-degree terms},
\end{equation*}
we get, for a monomial $P = \del^n / \del t_{k_1} \cdots \del t_{k_n}$,
\begin{equation*}
  D_P = c_{k_1,k_2,\dots,k_n} {s_1}^{\sum(2k_i+1)}
  {\del^n / \del s_{k_1} \cdots \del s_{k_n}} + \text{lower
    $s_1$-degree terms}.
\end{equation*}
So that, evaluating at $(0,1,0,0,\dots)$,
\begin{equation*}
  D_P(0,1,0,0,\dots;\del / \del{s_*}) = c_{k_1,k_2,\dots,k_n} 
  {\del^n / \del s_{k_1} \cdots \del s_{k_n}} + \text{lower
    order terms},
\end{equation*}
where the omitted terms are differential operators $\del_{s_*}^{m_*}$
such that $\norm{m_*}_+ < \norm{k_*}_+$.

Thus, in the bases \(\{\del^n/\del t_{k_1}\cdots\del t_{k_n}\}\) and
\(\{\del^m/\del s_{l_1}\cdots\del s_{l_m}\}\), the linear map \(P
\mapsto D_P\) is triangular and, therefore, invertible.
\end{proof}

\subsection{A matrix integral interpretation}
The Di~Francesco-Itzykson-Zuber theorem first appeared as a statement
about Hermitian matrix integrals; we recover the original formulation
by translating \prettyref{cor:dfiz-der} into the language of Gaussian
integrals related to the 't~Hooft-Kontsevich model.
\begin{corollary}[DFIZ Theorem] \label{cor:dfiz}
  There exists a vector space isomorphism 
\begin{equation*}
{Q}: \setC[\del / \del{t_*}] \to
  \setC [\tr X, \tr X^3, \tr X^5, \dots]
\end{equation*}
 such that, for $N \gg 0$,
  \begin{multline*}
    P(\partial_{t_*})\int_{\Hermitian[N]}\exp\Biggl\{\frac{\U}{6}
    \tr X^3\Biggr\}d\mu_{\Lambda}(X)=
    \int_{\Hermitian[N]}{Q}_P(X) \times 
    \\
    \times \exp\Biggl\{\frac{\U}{6}
    \tr X^3\Biggr\}d\mu_{\Lambda}(X)
  \end{multline*}
  in the sense of asymptotic expansions. 
\end{corollary}
The more general statement of \prettyref{cor:dfiz+der} corresponds to the
following.
\begin{corollary}
  \label{cor:dfiz+} 
  There exists a linear map 
  \begin{equation*}
    {Q^{s_*^\circ}}:\setC[\del / \del{t_*}] \to
    \setC[s_*^\circ; \tr X, \tr X^3, \ldots ]
  \end{equation*}
  such that, for $N\gg 0$,
  \begin{align*}
    P(\del / \del{t_*}) &\int_{\Hermitian[N]}
    \exp\Biggl\{-\U\sum_{j=0}^r (-1/2)^j s^\circ_j
    \frac{\tr(X^{2j+1})}{2j+1}\Biggr\}d\mu_{\Lambda}(X)=
    \\
    &=\int_{\Hermitian[N]}
    Q^{s_*^\circ}_P(s_*^\circ;X)\exp\Biggl\{-\U\sum_{j=0}^r (-1/2)^j
    s^\circ_j \frac{\tr (X^{2j+1})}{2j+1}\Biggr\}d\mu_{\Lambda}(X),
  \end{align*}
  in the sense of asymptotic expansions.
\end{corollary}

\subsection{A geometrical interpretation}

 As we have already remarked, differentiating the partition function
 $Z(s_*;t_*)$ with respect to the variable $t_k$ corresponds to
 ``inserting a $\tau_k$ in the coefficients''. Then, evaluating at the
 point $s_*^\circ=(0,1,0,0,\dots)$, one considers just the
 combinatorial class corresponding to ribbon graphs with only trivalent
 vertices, i.e., to the fundamental class in the moduli spaces
 $\overline{\mathcal M}_{g,n}$. This means that the action of
 $\text{ev}_{s_*^\circ}\circ \del/\del t_k$ on $Z(s_*;t_*)$ describes
 the linear functionals
 \begin{equation}
   \int_{\overline{\mathcal M}_{g,n}} {\psi_1}^k\land-\colon
   \setC[\psi_2,\dots,\psi_n]_{g,n}\to\setC,
 \end{equation}
 where $\setC[\psi_2,\dots,\psi_n]_{g,n}$ is the subalgebra of
 $H^*(\overline{\mathcal M}_{g,n})$ generated by the classes
 $\psi_2,\dots,\psi_n$.  More in general, the action of operators
 $\text{ev}_{s^\circ_*} \circ P(\del/\del t_*)$ on the partition
 function $Z(s_*;t_*)$ describes the linear functionals given by
 \begin{equation*}
   \int_{\overline{\mathcal M}_{g,n}} \tilde{P}(\psi_*)\land-
 \end{equation*}
 where $\tilde{P}(\psi_*)$ is a polynomial in the Miller classes.

 On the other hand, for any polyindex $m_*$, acting on $Z(s_*;t_*)$
 with the operator $\text{ev}_{s_*^\circ}\circ \del^{|m_*|}/{\del
   s_0}^{m_0}\cdots {\del s_r}^{m_r}$ corresponds to integrating the
 $\psi$ classes on the combinatorial stratum $W_{m_*;n}$ described by
 ribbon graphs having exactly $m_i$ distinguished $(2i+1)$-valent
 vertices, and all the other vertices of valence three. In other words,
 the action of $\text{ev}_{s_*^\circ}\circ \del^{|m_*|}/{\del
   s_0}^{m_0}\cdots {\del s_r}^{m_r}$ on $Z(s_*;t_*)$ describes the
 linear operators
 \begin{equation}
   \int_{W_{m_*;n}} \colon \setC[\psi_*]_{g,n}\to\setC.
 \end{equation}

 Therefore, \prettyref{cor:dfiz-der} could be interpreted by saying
 that the combinatorial classes and the $\psi$ classes define the same
 families of functionals on the subalgebra of the cohomology of the
 moduli spaces of stable curves generated by the Miller classes. So, in
 a certain sense (i.e., up to push-forwards and addition of classes
 supported in the boundary of moduli), one can say that the
 combinatorial classes are the Poincar{\'e} duals of the Mumford
 classes.  For a precise statement and more details on this topic, see
 \cite{kontsevich;feynman, arbarello-cornalba;dfiz,
   igusa;miller-morita, mondello;tautological}

\subsection{Example computation: $\del Z(t_*;s_*) / \del t_0$}
\label{sec:sample1}

Equation \eqref{eq:derivatives-partition2} tells us
\begin{equation}
  \label{eq:xmp1-start}
  \frac{\partial}{\partial t_0} \langle\!\langle\emptyset\rangle\!\rangle_{\Lambda, s_*} =
  -\sum_{\Gamma\in {\mathcal{S}_{1}}} \Coeff_{z}^{-1}
  \langle\!\langle\Gamma\rangle\!\rangle^{[1]}_{z\oplus\Lambda, s_*},
\end{equation}
where ${\mathcal{S}_{1}}$ is the set of hole types with a
$z$-decorated hole bounded only by one edge.  It consists of elements
\begin{equation*}
  {\xy,(0,-1)*{\bullet},(0,1)*{\bullet}%
    ,(0,-1);(0,1)**\dir{-},(0.3,0)*{z}\endxy}\quad  
  {\xy\rgvertex{11}\loose3\join7{10},(0,-.5)*{z}\endxy}
  {\xy\rgvertex{11}\loose3\loose2\loose4\join7{10},(0,-.5)*{z}\endxy}
  {\xy\rgvertex{11}\loose3\loose2\loose4\loose1\loose5\join7{10}%
    ,(0,-.5)*{z}\endxy}
  \quad\dots
\end{equation*}

The first graph in the list above has exactly \emph{two}
automorphisms, while none of the other graphs has non-trivial
automorphisms.  According to the Feynman rules, one computes
\begin{equation*}\label{eq:xmp1-1}
  \Coeff_z^{-1} \onehalf Z_{z\oplus\Lambda, s_*}
  \Biggl(\, {\xyc,(0,-1)*{\bullet},%
    (0,1)*{\bullet},(0,-1);(0,1)**\dir{-},(0.3,0)*{z}\endxyc} \Biggr)
  = -\onehalf{s_0}^2 = -\onehalf{s_0}^2 Z_{\Lambda}(\emptyset),
\end{equation*}
so that
\begin{equation*}
  \Coeff_z^{-1}
  \gint[z\oplus\Lambda,s_*]{\,\xyc,(0,-1)*{\bullet},%
    (0,1)*{\bullet}%
    ,(0,-1);(0,1)**\dir{-},(0.3,0)*{z}\endxyc}
  ^{[1]}
  = -\onehalf{s_0}^2 \langle\!\langle\emptyset\rangle\!\rangle_{\Lambda, s_*}^{}\,.
\end{equation*}
Moreover,
\begin{multline*}\label{eq:xmp1-2}
  \Coeff_z^{-1} 
  {\xyc%
    \rgvertex{11}%
    \fence1{i_1}{i_2\ }%
    \fence2{}{i_3}%
    \missing3%
    \fence4{i_{2m}}{i_{2m+1}}%
    \fence5{}{i_1}%
    \join7{10},%
    (0,-.5)*{z}%
    \endxyc} 
  =\Coeff_z^{-1} \left( -\U (-1/2)^{m+1}s_{m+1} \cdot \frac{2}{z+\Lambda_{i_1}} \right)
  \\
  = \frac{(-1)^m \U (2m+1)}{2^m}s_{m+1} \cdot  
  \left( \frac{1}{2m+1}
    {\xyc*!LC\xybox{\rgvertex[\!]{5}%
        \fence1{i_1}{i_2}%
        \fence2{i_2}{i_3}%
        \missing3%
        \fence4{i_{2m}}{\, i_{2m+1}}%
        \fence5{}{i_1}%
      }%
      \endxyc}
  \right),
\end{multline*}
so that 
\begin{multline*}
  \Coeff_z^{-1} \gint[z\oplus\Lambda,s_*]{  
    {\xyc%
      \rgvertex{11}%
      \loose1%
      \loose2%
      \missing3%
      \loose4%
      \loose5%
      \join7{10},%
      (0,-.5)*{z}%
      \endxyc} 
  }^{[1]}
  = \frac{(-1)^m \U (2m+1)}{2^m}s_{m+1} \cdot
  \langle\!\langle{\sf v}_{2m+1}\rangle\!\rangle_{\Lambda, s_*}^{}
  \\
  = -(2m+1)s_{m+1} \cdot \frac{\partial}{\partial s_m} \langle\!\langle \emptyset \rangle\!\rangle^{}_{\Lambda, s_*}
\end{multline*}

Therefore, equation \prettyref{eq:xmp1-start} becomes:
\begin{equation*}
  \label{eq:xmp1-final}
  \frac{\partial}{\partial t_0} 
  \langle\!\langle\emptyset\rangle\!\rangle_{\Lambda, s_*}^{} = 
  \Biggl( \frac{{s_0}^2}{2}+\sum_{m=0}^{\infty}(2m+1)s_{m+1}\frac{\del}{\del
    s_{m}}\Biggr) \langle\!\langle \emptyset
  \rangle\!\rangle^{}_{\Lambda, s_*};
\end{equation*}
we can rewrite it as:
\begin{equation*}
  \frac{\partial}{\partial t_0}
  Z(s_*;t_*)\biggr\rvert_{t_*(\Lambda)}=
  \Biggl(
  \frac{{s_0}^2}{2}+\sum_{m=0}^{\infty}(2m+1)s_{m+1}\frac{\del}{\del
    s_{m}}\Biggr)
  Z(s_*;t_*)\biggr\rvert_{t_*(\Lambda)}.
\end{equation*}
In the large $N$ limit, the $t_*(\Lambda)$ are independent
coordinates, thus
\begin{equation*}
  \frac{\partial}{\partial t_0}
  Z(s_*;t_*)=                 
  \Biggl(
  \frac{{s_0}^2}{2}+\sum_{m=0}^{\infty}(2m+1)s_{m+1}\frac{\del}{\del
    s_{m}}\Biggr) 
  Z(s_*;t_*),
\end{equation*}
which can be rewritten as:
\begin{equation}\label{eq:generates-relations}
  \frac{\partial}{\partial t_0}
  F(s_*;t_*)=  
  \frac{{s_0}^2}{2}+\sum_{m=0}^{\infty}(2m+1)s_{m+1}\frac{\del}{\del
    s_{m}}
  F(s_*;t_*).
\end{equation}
Now, by equation \prettyref{eq:generates-relations}, one finds:
\begin{equation*}
  \frac{\partial^3}{\partial {t_0}^3}
  F(t_*)=1+\left.\left(\frac{\del^3}{\del
        s_{0}^3}
      F(s_*;t_*)\right)\right\rvert_{s_*=(0,1,0,\dots)},
\end{equation*}
which implies
\begin{equation*}
  \langle\tau_0^3\tau_{\nu_1}\cdots\tau_{\nu_n}\rangle=
  1+\sum_{m_1=0}^\infty\langle\tau_{\nu_1}\cdots\tau_{\nu_n}\rangle_{(3,m_1,
    0,\dots);n}.
\end{equation*}
In particular, for $n=0$, one recovers the well-known relation
$\langle{\tau_0}^3\rangle=1$
\cite{kontsevich;intersection-theory;1992,
  witten;2dgravity}.

\subsection{Example computation: $\del
  \langle\!\langle\emptyset\rangle\!\rangle_{\Lambda, s_*}^{} / \del t_1$ at
  $s_* = (0,0,s_2, 0,\ldots)$}
\label{sec:sample2}

As an illustration of \prettyref{cor:dfiz+}, we compute $\del
\langle\!\langle\emptyset\rangle\!\rangle_{\Lambda, s_*}^{} / \del
t_1$ at $s^\circ_* = (0, 0, s_2, \break 0, \ldots)$.  Since \(s_i=0\) for
\(i\neq 2\) we need to consider graphs with 5-valent vertices only;
moreover, since $k = 1$, we need to consider only holes made up of at
most $3$ edges. The relevant hole types therefore are:
\begin{align*}
  &\Gamma_{1} := {\xyc
    \rgvertex{11}%
    \join{7}{10}\loose2\loose3\loose4%
    ,(0,-.5)*{\scriptstyle{z}}
    \endxyc}
  &&
  \Gamma_{2} := {\xyc
    \rgvertex{6}\loose4\loose5\loose3
    ,(0,0);(2,0)**\crv{(1,1)}
    ,(0,0);(2,0)**\crv{(1,-1)}
    ,(2,0)\rgvertex{6}\loose1\loose2\loose6,(1,0)*{\scriptstyle{z}}
    \endxyc}
  \\
  &\Gamma_{3} := {\xyc
    \rghole[z]{3}\loose1\loose2\loose3
    ,(0,1);(.5,1.4)**\dir{-}
    ,(0,1);(-.5,1.4)**\dir{-}
    ,(-.866,-.5);(-.98,-1)**\dir{-}
    ,(-.866,-.5);(-1.32,-.28)**\dir{-}
    ,(.866,-.5);(.98,-1)**\dir{-}
    ,(.866,-.5);(1.32,-.28)**\dir{-}
    \endxyc}
  &&
  \Gamma_{4} := {\xyc
    \rgvertex{6}\loose3\loose5\join{2}{6}
    ,(3,0)\rgvertex{6}\loose2\loose6\join{3}{5}
    ,(0,0);(.6,0)**\dir{-}
    ,(.8,0);(2.2,0)**\dir{-}
    ,(2.4,0);(3,0)**\dir{-}
    ,(1.5,.2)*{\scriptstyle{z}}
    \endxyc}
\end{align*}
By equation \prettyref{eq:dZst-twotwo},
\begin{equation}\label{eq:eqesempio}
  \left.
    \frac{\del
      \langle\!\langle\emptyset\rangle\!\rangle_{\Lambda, s_*}^{}}{\del
      t_1}
  \right\vert_{%
    s_*=(0,0,s_2,0,\dots) 
  }
  =\left. -\sum_{i=1}^4 \Coeff_z^{-3}\langle\!\langle
    \Gamma_{i}
    \rangle\!\rangle_{z\oplus\Lambda, s_*}^{[1]}\right\vert_{%
    s_*=(0,0,s_2,0,\dots) 
  }.
\end{equation}

One computes:
\begin{equation*}
  \Coeff_z^{-3}\left(
    {\xyc
      \rgvertex{11}%
      \join{7}{10}%
      \fence2{i}{j}%
      \fence3{j}{k}%
      \fence4{\!\!k}{i}%
      ,(0,-.5)*{\scriptstyle{z}}
      \endxyc}
  \right)
  = - \frac{\U}{2} s_2\,{\Lambda_i}^2
  = -\frac{\U}{2} s_2  
  {\xyc%
    \rgvertex[\varphi_1]{3}%
    \fence1{i}{j}
    \fence2{j}{k}
    \fence3{k}{i}
    ,(.5,.4)*{\scriptstyle{\bigstar}}%
    \endxyc}
  ,
\end{equation*}
with \(\varphi_1 (\theta_1) = {\theta_1}^2\), so that
\begin{equation*}
  \Coeff_z^{-3}\langle\!\langle\Gamma_{1}\rangle\!
\rangle^{[1]}_{z\oplus\Lambda, s_*}
  = -\frac{\U}{2} s_2
  \gint{ 
    {\xyc%
      \rgvertex[\varphi_1]{3}%
      \loose1%
      \loose2%
      \loose3%
      ,(.5,.4)*{\scriptstyle{\bigstar}}%
      \endxyc}
  }.
\end{equation*}
The polynomial $\varphi_1$ is not cyclically invariant, but we can change
it into a cyclically invariant one by means of formula
\prettyref{eq:to-cyclic}:
\begin{equation*}
  \gint{%
    {\xyc%
      \rgvertex[\varphi_1]{3}%
      \loose1%
      \loose2%
      \loose3%
      \endxyc}
  }
  = \gint{
    {\xyc%
      \rgvertex[\overline{\varphi}_1]{3}%
      \loose1%
      \loose2%
      \loose3%
      \endxyc}},
  \qquad
  \overline{\varphi}_1(\theta_1,\theta_2,\theta_3) =
  {\theta_1}^2 + {\theta_2}^2 +
  {\theta_3}^2.
\end{equation*}
To avoid cumbersome notations, let us write 
\begin{equation*}
  \langle\!\langle\Gamma\rangle\!\rangle_{s_*^\circ}:=\langle\!\langle\Gamma\rangle\!\rangle_{\Lambda, s_*}^{}\Bigr\rvert_{s_*=(0,
    0,
    s_2, 0, \ldots)}\,,
\end{equation*}
so that, evaluating at \(s_*=(0, 0, s_2, 0, \ldots)\), we find:
\begin{equation*}
  \label{eq:cumbersome}
  \Coeff_z^{-3}\langle\!\langle\Gamma_{1}\rangle\!
  \rangle^{[1]}_{z\oplus\Lambda, s_*}\Bigr\rvert_{s_*=(0, 0,
    s_2, 0, \ldots)}
  = -\frac{\U}{2} s_2
  \langle\!\langle {\sf v}_3^{\overline{\varphi}_1}
  \rangle\!\rangle_{s_*^\circ}.
\end{equation*}
In the same way (and accounting for the
automorphism groups involved) we get:
\begin{align*}
  \Coeff_z^{-3}\langle\!\langle\Gamma_{2}\rangle\!
  \rangle^{[1]}_{z\oplus\Lambda, s_*}\Bigr\rvert_{s_*=(0, 0,
    s_2, 0, \ldots)}
  &= \frac{{s_2}^2}{4} \langle\!\langle
  {\sf
    v}^{\overline{\varphi}_2}_{6}\rangle\!\rangle_{s_*^\circ}\,,
  \intertext{where
    $\overline{\varphi}_2(\theta_1,\theta_2,\dots,\theta_6) =
    {\theta_1} +{\theta_2}+\dots+{\theta_6}$;}
  \Coeff_z^{-3}\langle\!\langle\Gamma_{3}\rangle\!
  \rangle^{[1]}_{z\oplus\Lambda, s_*}\Bigr\rvert_{s_*=(0, 0,
    s_2, 0, \ldots)}
  & = \frac{3\U}{8} {s_2}^3 \langle\!\langle
  {\sf v}^{}_{9}\rangle\!\rangle_{s_*^\circ}\,;
  \\
  \Coeff_z^{-3}\langle\!\langle\Gamma_{4}\rangle\!
  \rangle^{[1]}_{z\oplus\Lambda, s_*}\Bigr\rvert_{s_*=(0, 0,
    s_2, 0, \ldots)}
  & = -{s_2}^2 \langle\!\langle
  {\sf v}^{}_{2}{\textstyle\coprod} {\sf
    v}^{}_{2}\rangle\!\rangle_{s_*^\circ}\,.
\end{align*}
Now, equation \prettyref{eq:eqesempio} can be rewritten as:
\begin{multline}\label{eq:eqesempio2}
  \left.
    \frac{\del
      \langle\!\langle\emptyset\rangle\!\rangle_{\Lambda, s_*}^{}}{\del
      t_1}
  \right\vert_{%
    s_*=(0,0,s_2,0,\dots) 
  }
  =\frac{\U}{2} s_2 \langle\!\langle 
  {\sf v}_3^{\overline{\varphi}_1}
  \rangle\!\rangle_{s_*^\circ} 
  -\frac{{s_2}^2}{4} \langle\!\langle
  {\sf
    v}^{\overline{\varphi}_2}_{6}\rangle\!\rangle_{s_*^\circ}
  \\
  -\frac{3\U}{8}{{s_2}^3}
  \langle\!\langle
  {\sf v}^{}_{9}\rangle\!\rangle_{s_*^\circ}
  + {s_2}^2\langle\!\langle
  {\sf v}^{}_{2}{\textstyle\coprod} {\sf
    v}^{}_{2}\rangle\!\rangle_{s_*^\circ}\,.
\end{multline}

According to the proof of \prettyref{theorem:main}, we could forget the
contribution coming from the last term in the right-hand side, because
it contains even-valent residual vertices. However, we will not do
this, so to explicitly show how it gets canceled out.

Let us proceed to lower the degree of the polynomials decorating the
vertices in equation \prettyref{eq:eqesempio2} by contraction of
edges, starting with the trivalent vertex decorated by
\(\overline{\varphi}_1\). It is cyclically decomposable; a possible
decomposition is:
\begin{equation*}
  \overline{\varphi}_1(\theta_1, \theta_2, \theta_3) =
  (\theta_1 +
  \theta_2) \cdot
  \psi_1(\theta_1, \theta_2, \theta_3) + \text{cyclic permutations},
\end{equation*}
where
\begin{equation*}
  \psi_1(\theta_1, \theta_2, \theta_3) = \frac {\theta_1 + \theta_2 -
    \theta_3} {2}.
\end{equation*}
As in the proof of \prettyref{lemma:lower-degree},
\begin{align*}
\langle\!\langle{\sf v}_3^{\bar
\varphi_1}\rangle\!\rangle_{s^\circ_*}
  &= \gint[s^\circ_*]{{\xyc\rgvertex[\overline{\varphi}_1]{3}%
      \loose1\loose2\loose3\endxyc}}
  = \gint[s^\circ_*]{{\xyc\rgvertex[u_{\psi_1}]{3}%
      \loose1\loose2\loose3\cilialontano{1}{3}%
      \endxyc}}
  \\
  &= 
  \gint[s^\circ_*]{{\xyc\rgvertex[u_{\psi_1}]{3}%
      \loose1\loose2\loose3\cilialontano{1}{3}%
      ,(1.12,-0.64);(1.54,-1.07)**\dir{-}%
      ,(1.12,-0.64);(1.5,-.3)**\dir{-}%
      ,(1.12,-0.64);(1.7,-.6)**\dir{-}%
      ,(1.12,-0.64);(1.2,-1.2)**\dir{-}%
      ,(1.12,-0.64)*{\bullet}
      \endxyc}}
  +\gint[s^\circ_*]{{\xyc\rgvertex[u_{\psi_1}]{3}%
      \joinlontano{1}{3}\loose2\cilialontano{1}{3}\endxyc}}
  +\gint[s^\circ_*]{{\xyc\rgvertex[u_{\psi_1}]{3}%
      \joinlontano{3}{2}\loose1\cilialontano{1}{3}\endxyc}}
  \\
  &= 
  -\frac{\U}{2}s_2\gint[s^\circ_*]{{\xyc\rgvertex[\psi_1]{6}\cilia12%
      \loose1\loose2\loose3\loose4\loose5\loose6\endxyc}}
  + \tr\Lambda \gint[s^\circ_*]{\negquad 
    {\xyc\rgvertex[\!]{3}%
      \loose1\endxyc} 
    \negquad }
  \\
  &\hskip5em\relax
  - \tr\Lambda \gint[s^\circ_*]{\negquad
    {\xyc\rgvertex[\!]{3}%
      \loose1\endxyc}
    \negquad }
  + 2 \tr\Lambda^0 \gint[s^\circ_*]{\negquad
    {\xyc\rgvertex[\psi_2]{3}%
      \loose1\endxyc}
    \negquad},
\end{align*}
with $\psi_2(\theta_1) = \theta_1$. Find a cyclic equivalent of
\(\psi_1\), by applying \prettyref{eq:to-cyclic}
again:
\begin{equation*}
  \gint[s^\circ_*]{{\xyc\rgvertex[\psi_1]{6}\cilia12%
      \loose1\loose2\loose3\loose4\loose5\loose6\endxyc}}
  = \gint[s^\circ_*]{
    {\xyc\rgvertex[\overline{\psi}_1]{6}%
      \loose1\loose2\loose3\loose4\loose5\loose6\endxyc}
  }.
\end{equation*}
Explicit computation shows that $\overline{\psi}_1 = (1/2)
\overline{\varphi}_2$, so that:
\begin{equation*}
  \langle\!\langle {\sf
    v}^{\overline{\varphi}_1}_3\rangle\!\rangle_{s_*^\circ}
  =   - \frac{\U}{4} s_2
  \langle\!\langle {\sf
    v}^{\overline{\varphi}_2}_6\rangle\!\rangle_{s_*^\circ}
  + 2\tr(\Lambda^0)
  \langle\!\langle {\sf
    v}^{\psi_2}_1\rangle\!\rangle_{s_*^\circ},
\end{equation*} 
and both vertices on the right are decorated by cyclic polynomials.
Since \(\psi_2(\theta_1) = (\theta_1 + \theta_1) \cdot (1/2) =
u_{1/2}\), then \prettyref{eq:C1} gives:
\begin{equation}\label{eq:esempio2-two}
  \langle\!\langle {\sf
    v}^{\psi_2}_1\rangle\!\rangle_{s_*^\circ}
  = -\frac{\U}{4} {s_2} \gint[s^\circ_*]{
    {\xyc\rgvertex[\!]{4}\loose1\loose2%
      \loose3\loose4\cilia{1}{2}\endxyc}}
  = -\U {s_2} \langle\!\langle {\sf
    v}^{}_4\rangle\!\rangle_{s_*^\circ}\,,
\end{equation}
which is residual.

The other positive degree polynomial appearing on the right hand side
of equations \prettyref{eq:eqesempio2} and \prettyref{eq:esempio2-two}
is \(\overline{\varphi}_2\), which has a cyclic decomposition
\begin{equation*}
  \overline{\varphi}_2(\theta_1, \dots, \theta_6) =
  u_{1/2} +
  \text{cyclic
    permutations}.
\end{equation*} 
Again as in the proof of \prettyref{lemma:lower-degree}, we get:
\begin{align*}
  \langle\!\langle {\sf
    v}^{\overline{\varphi}_2}_6\rangle\!\rangle_{s_*^\circ}
  &=
  \gint[s^\circ_*]{{\xyc\rgvertex[u_{\frac{1}{2}}]{6}%
      \loose1\loose2\loose3\loose4\loose5\loose6\cilia{1}{2}\endxyc}}
  \\
  &= 
  \gint[s^\circ_*]{{\xyc\rgvertex[u_{\frac{1}{2}}]{6}%
      \loose1\loose2\loose3\loose4\loose5\loose6\cilia{1}{2}
      ,(1.3,0);(1.6,.5)**\dir{-},(1.3,0);(1.6,-.5)**\dir{-}%
      ,(1.3,0);(1.8,.2)**\dir{-},(1.3,0);(1.8,-.2)**\dir{-}%
      ,(1.3,0)*{\bullet}\endxyc}}+
  \gint[s^\circ_*]{{\xyc\rgvertex[u_{\frac{1}{2}}]{6}%
      \join{1}{2}\loose3\loose4\loose5\loose6\cilia{1}{2}\endxyc}}
  +\gint[s^\circ_*]{{\xyc\rgvertex[u_{\frac{1}{2}}]{6}%
      \joinunder{1}{2}{3}\loose2\loose4\loose5\loose6\cilia{1}{2}\endxyc}}
  \\
  &\qquad
  +\gint[s^\circ_*]{{\xyc\rgvertex[u_{\frac{1}{2}}]{6}\loose2\loose3\loose5%
      \loose6\cilia{1}{2}%
      ,(0.4,0);(.52,.8)**\crv{(1.3,0)},(-0.4,0);(-.52,.8)**\crv{(-1.3,0)}%
      ,(.44,.85);(-.44,.85)**\crv{(0,1.25)}\endxyc}}
  +\gint[s^\circ_*]{{\xyc\rgvertex[u_{\frac{1}{2}}]{6}%
      \joinunder{1}{6}{5}\loose2\loose3\loose4\loose6\cilia{1}{2}\endxyc}}
  +\gint[s^\circ_*]{{\xyc\rgvertex[u_{\frac{1}{2}}]{6}%
      \join{1}{6}\loose2\loose3\loose4\loose5\cilia{1}{2}\endxyc}}
  \\
  &= -\frac{\U}{4}s_2\gint[s^\circ_*]{{\xyc\rgvertex[\!]{9}%
      \loose1\loose2\loose3\loose4\loose5\loose6\loose7\loose8\loose9%
      \cilia{1}{2}\endxyc}} 
  +2\cdot\gint[s^\circ_*]{
    {\xyc\rgvertex[\!]{3}%
      \loose1\loose2\loose3\ciliamedio{1}{3}\endxyc}
    \negqquad
    {\xyc\rgvertex[\!]{3}%
      \loose1\endxyc}
    \negquad\negquad}
  \\
  &\qquad
  + 2\cdot\gint[s^\circ_*]{\negquad 
    {\xyc\rgvertex[\!]{6}%
      \loose3\loose5\cilia{1}{2}\endxyc}
    \negqquad\negqquad
    {\xyc\rgvertex[\!]{6}%
      \loose2\loose6\cilia{1}{2}\endxyc}
    \negquad}
  + 2\tr\Lambda^0 \gint[s^\circ_*]{{\xyc\rgvertex[\!]{4}%
      \loose1\loose2\loose3\loose4\cilia{1}{2}\endxyc}}
  \\
  &= -\frac{9\U}{4}s_2\langle\!\langle {\sf
    v}^{}_9\rangle\!\rangle_{s_*^\circ}
  + 8\cdot\tr\Lambda^0 \langle\!\langle {\sf
    v}^{}_4\rangle\!\rangle_{s_*^\circ}
  + 6\,\langle\!\langle {\sf
    v}^{}_3{\textstyle\coprod} {\sf v}^{}_1\rangle\!\rangle_{s_*^\circ}
  + 8\,\langle\!\langle {\sf
    v}^{}_2 {\textstyle\coprod} {\sf v}^{}_2\rangle\!\rangle_{s_*^\circ} .
\end{align*}

In the end, we substitute back into \prettyref{eq:eqesempio2}:
\begin{equation*}
  \left.
    \frac{\del
      \langle\!\langle\emptyset\rangle\!\rangle_{\Lambda, s_*}^{}}{\del
      t_1}
  \right\vert_{%
    s_*=(0,0,s_2,0,\dots) 
  }
  \negquad= -\frac{3}{4} {s_2}^2  \langle\!\langle {\sf
    v}^{}_3 {\textstyle\coprod} {\sf v}^{}_1\rangle\!\rangle_{s_*^\circ}
  -\frac {3\U} {32} {s_2}^3 \langle\!\langle {\sf
    v}^{}_9\rangle\!\rangle_{s_*^\circ}.
\end{equation*}
In terms of Hermitian matrix integrals this reads:
\begin{multline*}
  \frac{\del}{\del t_1}\int_{\Hermitian[N]} \exp\Biggl\{ -\frac{\U}{4}
  s_2 \frac{\tr X^5}{5} \Biggr\} \ud\mu_{\Lambda}(X) =
  \int_{\Hermitian[N]} \Biggl( -\frac{1}{4} {s_2}^2 \tr X^3 \tr X +
  \\
  -\frac{\U}{96} {s_2}^3 \tr X^9 \Biggr) \exp\Biggl\{ -\frac{\U}{4}
  s_2 \frac{\tr X^5}{5} \Biggr\} \ud\mu_{\Lambda}(X).
\end{multline*}

\begin{acknowledgement}
  The authors thank their thesis advisor Enrico Arbarello, who has
  constantly been of assistance and support during the (long)
  preparation of this paper, and has had the unrewarding task of
  reading the unreadable early drafts. We owe special thanks also to
  Gilberto Bini, Federico De~Vita, and Gabriele Mondello, for their
  interest and continuous encouragement.  Constructive criticism by
  the Referees helped us improve the exposition.
  
  Last, we thank Kristoffer H. Rose and Ross Moore, the authors of
  \Xy-Pic: without this wonderful \TeX\ graphics package this paper
  would probably have never been written.
\end{acknowledgement}


\end{document}